\magnification  1200
\ifx\eplain\undefined \input eplain \fi
\baselineskip13pt 

\overfullrule=0pt


\def\wtilde{\widetilde}

\def\bar{\overline}
\font\smalsmalbf=cmbx8



\font\smalltenrm=cmr8
\font\smallteni=cmmi8
\font\smalltensy=cmsy8
\font\smallsevrm=cmr6   \font\smallfivrm=cmr5
\font\smallsevi=cmmi6   \font\smallfivi=cmmi5
\font\smallsevsy=cmsy6  \font\smallfivsy=cmsy5
\font\smallsl=cmsl8      \font\smallit=cmti8

\def\smallfonts{\lineadj{80}\textfont0=\smalltenrm  \scriptfont0=\smallsevrm
                \scriptscriptfont0=\smallfivrm
    \textfont1=\smallteni  \scriptfont1=\smallsevi
                \scriptscriptfont0=\smallfivi
     \textfont2=\smalltensy  \scriptfont2=\smallsevsy
                \scriptscriptfont2=\smallfivsy
      \let\it\smallit\let\sl\smallsl\smalltenrm}


\font\mathbold=msbm9 at 10pt

\def \C{{\hbox{\mathbold\char67}}}
\def\N{{\hbox{\mathbold\char78}}}
\def\R{{\hbox{\mathbold\char82}}}

\def\e{\epsilon}

\def\un{{1\over n-2}}

\def\imp{{\;\;\Longrightarrow\;\;}}
\def\lineadj#1{\normalbaselines\multiply\lineskip#1\divide\lineskip100
\multiply\baselineskip#1\divide\baselineskip100
\multiply\lineskiplimit#1\divide\lineskiplimit100}

\def\remark#1.{\medskip{\noin\bf Remark #1.\enspace}}
\def\endpf{$$\eqno/\!/\!/$$}

\def\pf#1.{\smallskip\noin{\bf  #1.\enspace}}

\def\noin{\noindent}

\def\ds{\displaystyle}
\def\ts{\textstyle}

\def\e{\epsilon}\def\part{\partial_t}

\def\wtilde{\widetilde}

\def\Rn{\R^n}\def\B{{\cal B}}

\def\ref#1{{\bf{[#1]}}}

\def\p{\partial}

\def\half{{\ts{1\over2}}}

\def\H{{\cal H}}

\centerline{\bf  Sharp  Moser-Trudinger inequalities for the Laplacian without boundary conditions} 
\bigskip
\centerline {Luigi Fontana, Carlo Morpurgo}
\vskip1em
\midinsert
\noindent{\smalsmalbf Abstract. }{\smallfonts We derive a sharp Moser-Trudinger inequality for the borderline Sobolev imbedding
of $W^{2,n/2}(B_n)$ into the exponential class, where $B_n$ is the unit ball of $\R^n$. The corresponding sharp results for the spaces $W_0^{d,n/d}(\Omega)$ are well known, for general domains $\Omega$, and are due to Moser and Adams. When the zero boundary condition is removed the only known results are for $d=1$ and are due to Chang-Yang, Cianchi and Leckband. Our proof is based on general abstract results recently obtained by the authors in [FM], and on a new integral representation formula for the ``canonical" solution of the Poisson equation on the ball, that is the unique solution of the equation $\Delta u=f$ which is orthogonal to the harmonic functions on the ball. The main technical difficulty of the paper is to establish an asymptotically sharp growth estimate for  the kernel of such representation, expressed in terms of its distribution function. We will also consider the  situation where the exponential class is endowed with 
more general Borel measures, and obtain corresponding sharp  Moser-Trudinger inequalities of trace type.

}
\def\bar{\overline}
\endinsert\bigskip\vskip2em
\centerline{\bf 1. Prologue}\bigskip\bigskip

A Moser-Trudinger inequality is  a statement about the exponential integrability of functions belonging to  the Sobolev space $W^{k,n/k}(\Omega)$, where $\Omega$ is an open set of an $n-$dimen- sional manifold, and $1\le k<n$. In general terms, suppose that $\nu$ is a  Borel measure on $\Omega$ with $\nu(\Omega)<\infty$ and $P_k$ is a  differential (or pseudodifferential) operator of order $k$, acting on a subspace $H$ of $W^{k,n/k}(\Omega)$ so that $P_ku\neq 0$ if $u\in H$, unless $u=0$. In this situation, establishing a  sharp Moser-Trudinger inequality, in its basic form, consists in proving the existence of an optimal  constant $\alpha>0$ for which 
$$\sup_{u\in H}\;\int_\Omega \exp\bigg[\alpha\bigg(
{|u(x)|\over\|P_k u\|_{n/k}}\bigg)^{n\over n-k}\bigg]
\,d\nu(x)<\infty.\eqno(0)$$
A wealth of results exist for $H=W_0^{k,n/k}(\Omega)$, the closure of $C_0^\infty (\Omega)$ in $W^{k,n/k}(\Omega)$, or  when $\Omega$  is itself a compact manifold without boundary, in which case obviously  $W_0^{k,n/k}(\Omega)=W^{k,n/k}(\Omega)$. In the case of  bounded $\Omega\subseteq \Rn$, endowed with the Lebesgue measure, the first sharp result is due to Moser [Mo], for  $k=1$ and $P_1=\nabla$, the classical gradient operator. This result was later extended by Adams in [Ad], to integer powers of the Laplacian and their gradients; many more extensions, generalizations and variations of Adams' and Moser's results have appeared since (for a partial list see for example the cited works in  [Ci1], [Ci2] and [FM]). The present authors recently unified and improved Adams' strategy to  a general measure-theoretic setting, and provided several new  sharp inequalities of type (0), for rather general operators $P_k$ and measures $\nu$, with $H=W_0^{k,n/k}(\Omega)$, both on Riemannian and subRiemannian manifolds ([BFM], [FM]).
\eject
In contrast, not much is known  about inequality (0) for functions $u\in W^{k,n/k}(\Omega)$ that do not necessarily vanish on the boundary of $\Omega$, that is, when $H$ is allowed to contain functions that do not necessarily belong to $W_0^{k,n/k}(\Omega)$. So far the only results available are for the case $k=1$, when $P_k=\nabla$, on a certain class of domains in $\Rn$. In this situation the obvious candidate for $H$ is the space of functions of $W^{1,n}(\Omega)$ with zero mean, that is functions orthogonal to the constants. The most general result can be stated roughly as follows: suppose that $\Omega\subseteq \Rn$ is a bounded domain of class $C^{1,\alpha}$ except for finitely many conical singularities at the boundary; let $\theta_\Omega$ be the minimum aperture of the  cones at those  singularities. Then, there exists a constant $C$ such that  for all $u\in W^{1,n}(\Omega)$ (except of course for the constant  function 0)
$$\int_{\Omega}\exp\bigg[{n(\theta_\Omega)^{1\over {n-1}}}\bigg({|u(x)-u_\Omega|\over\|\nabla u\|_n}\bigg)^{n\over n-1}\bigg]dx\le C,\eqno(1) $$
where  $u_\Omega$ is the average of $u$ over $\Omega$.
In case of $C^{1,\alpha}$ (in particular smooth) domains
 we clearly  have $\theta_{\Omega}={\half}\omega_{n-1}$ where 
$\omega_{n-1}$
is the surface measure of the unit sphere of $\Rn$. The first version of  this result is due to Chang and Yang [CY], and dates back to 1988, for  piecewise $C^2$ domains of $\R^2$. The $n-$dimensional extension given above  was found by Cianchi [Ci1] in 2005, and independently by Leckband [Le], but only for the unit ball of $\Rn$. 

It is not difficult to realize that the sharp constant in (1) has to be smaller than $n(\theta_\Omega)^{1\over {n-1}}$; this was already observed by Fontana in 1993 [Fo]. The classical sharp Moser-Trudinger inequality for $W_0^{1,n}(\Omega)$ is extremized by a family of functions $u_r$, the so-called Moser functions, which are radial and centered at an interior point. This means that $u_r\in W_0^{1,n}(\Omega)$, and the functional in (1) along this family
can be made arbitrarily large  if the exponential constant is greater than $n\omega_{n-1}^{1\over n-1}$, as $r\to 0$. On the other hand, if $\wtilde u_r$ denote the same functions but centered at a boundary point,
then $\wtilde u_r\in W^{1,n}(\Omega)$ and it is not hard to check that $\|\nabla \wtilde u_r\|_n^n\sim (\theta_\Omega/\omega_{n-1})\|\nabla u_r\|_n^n$, as $r\to 0$, whence the family $\wtilde u_r$ extremizes (1).
 
More recently,  Cianchi further extended (1) to a general class of  Borel measures $\nu$, obtaining trace-type inequalities,  by allowing  $u$ to belong to more general Lorentz-Sobolev spaces, and using regularizing functions $u_\Omega$ more general than the average. In another direction,  Pankka, Poggi-Corradini, and Rajala [PPCR] derive a sharp trace version of (1) on the boundary of the unit ball $B_n$ of $\Rn$, where the functions $u$ involved are those in   $W^{1,n}(B_n)$  that are continuous, monotone and with $u(0)=0$.

It is natural to speculate that there should be sharp  versions  of (1) for operators of order higher than 1, however at present there are no published results of this sort, not even in the simplest Euclidean settings.
The purpose of the present paper is to  give a complete answer to this problem for the simplest operator of order 2, the Laplacian, on the simplest smooth Euclidean domain, the unit ball. 
\bigskip\medskip\eject

\centerline{\bf 2. Statements of main results}\bigskip
 Let us set some notation.  Let $B_n=\{x\in\Rn:|x|<1\}$ denote the open unit ball of $\Rn$  and $S^{n-1}=\partial B_n$  the unit sphere;   denote their volumes in the corresponding standard Euclidean metrics by 
$$|S^{n-1}|=\omega_{n-1}={2\pi^{n/2}\over\Gamma\big({n\over2}\big)},\qquad |B_n|={\omega_{n-1}\over n}.$$ 
We will also denote the open ball of center $a$ and radius $r$ by $B(a,r)$.\medskip
The usual Sobolev space on an open set $\Omega$ is denoted as $W^{k,p}(\Omega)$, the set of functions in $L^p(\Omega)$ whose distributional derivatives $D^\gamma u$, up to order $k$ are also in 
$L^p(\Omega)$. A norm in $W^{k,p}(\Omega)$ is given as $\|u\|_{k,p}=\Big(\sum_{|\gamma|\le k} \|D^\gamma u\|_p^p\Big)^{1/p}$. The space $W_0^{k,p}(\Omega)$ is the closure of $C_0^\infty(\Omega)$ in $\big(W^{k,p}(\Omega),\|\cdot\|_{k,p}\big)$.

The standard Laplacian is the operator $\Delta=\sum_1^n\partial_{jj}^2$ and its fundamental solution
for $n\ge 3$ is given by the Newtonian kernel 
$$N(x)=-c_n |x|^{2-n},\qquad
c_n={1\over(n-2)\omega_{n-1}}\eqno(2)$$
and for $n=2$
$$N(x)={1\over 2\pi}\log|x|,\eqno(3)$$
in the sense that $\Delta N(\cdot-y)=\delta_y$, the Dirac delta at $y$.

To better describe our results let us first recall a special case of the Adams sharp inequality for the Laplacian: for any open and bounded $\Omega$ in $\Rn$,$\;n\ge3$,  there exists $C>0$ such that for all $u\in W_0^{2,n/2}(\Omega)$
$$\int_{\Omega}\exp\bigg[{c_n^{-{n\over n-2}}\over |B_n|}\bigg({|u(x)|\over \|\Delta u\|_{n/2}}\bigg)^{n\over n-2}\bigg]dx\le C\eqno(4)$$
where the exponential constant in (4) is sharp, i.e. it cannot be replaced by a larger constant.

 Our goal is to establish a version of (4) when the zero boundary condition is removed and when $\Omega=B_n$, that is for functions in $W^{2,n/2}(B_n)$. Clearly this imposes some restrictions on  the function $u$, which  is not allowed to be harmonic. In analogy with the case $k=1$ it is natural to impose the condition that our functions $u$ be orthogonal to the space of $L^2$ harmonic functions, the so-called $L^2$ Bergman space. For $n\ge3$ this is actually possible, since by the classical embedding theorem $W^{2,n/2}$ is in any $L^q$, $\, n/2\le q<\infty$, and  hence it is in $L^2$. 
\eject
To be more specific, the $L^p$ harmonic Bergman space on $B_n$ is defined for $p\ge 1$ as
$$b^p=\{u\in L^p(B_n): \;\; u {\hbox{ harmonic in }} B_n\,\}$$
which is a closed subspace $L^p(B_n)$. In particular  $b^2$ is a closed subspace of $L^2(B_n)$ with the usual inner product $\langle u,v\rangle=\int_{B_n} u\bar v$. The harmonic Bergman projection is the unique orthogonal projection $R:L^2(B_n)\to b^2$, with kernel $R(x,y)$, the Bergman Kernel. One can show that the operator $R$ with kernel $R(x,y)$ can in fact be defined on any $L^p$, and $R:L^p(B_n)\to b^p$, for $1<p<\infty$. Moreover $R:W^{k,p}(B_n)\to W^{k,p}(B_n)$, if $p>1$,  and it is bounded (see e.g. [KK], Thm. 4.5). 
The $L^2$-orthogonal of the space $b^2$ will be denoted as $(b^2)^{\perp}$, and it's clear that since $R:W^{2,n/2}\to W^{2,n/2}$ then 
$$(b^2)^\perp\cap W^{2,n/2}=\big\{u-Ru,\;u\in W^{2,n/2}\big\}.$$  
More generally, if $p>1$ and $p'$ is the conjugate to $p$ define the subspace of $L^p(B_n)$
$$\big(b^{p'}\big)^\perp=\big\{v\in L^p(B_n):\;\int_{B_n} v\bar h=0,\, \forall h\in b^{p'}\big\}$$
and we have, by density, that   $(b^{p'})^\perp\cap W^{2,p}=\big\{u-Ru,\;u\in W^{2,p}\big\}.$

These basic  facts are true on any smooth domain, however on the ball the Bergman kernel can be explicitly computed, and  this is indeed the main reason why we work on the unit ball.

We are now ready to state our main theorem.\smallskip
\proclaim Theorem 1. For any $n\ge 3$ there exists a constant $C$ such that for any $u\in W^{2,n/2}(B_n)$ and $u\notin b^2$
$$\int_{B_n}\exp\bigg[\alpha_n\bigg({|u(x)-Ru(x)|\over \|\Delta u\|_{n/2}^{}}\bigg)^{n\over n-2}\bigg]dx\le C\eqno(5)$$
where 
$$\alpha_n=\cases {\ds{1\over |B_n|}\,c_n^{-{n\over n-2}} & if $n=3$ or $n=4$\cr 
\ds{{1\over |G_n|}}\,c_n^{-{n\over n-2}}& if $n\ge 5$\cr}\eqno(6)$$
where 
$$|G_n|={\pi^{n-1\over2}\over n}(2n-4)^{n\over n-2}{\Gamma\big({1\over2}+{n\over n-2}\big)\over\Gamma\big({n\over2}+{n\over n-2}\big)}\eqno(7)$$
is the volume of the $n$-dimensional convex region
$$G_n=\{y=(y_1,...,y_n)\in\Rn:\;y_1>0,\,|y|^n< (2n-4)y_1^2\}.$$
The constant $\alpha_n$ in (5)  is sharp, in the sense that it cannot be replaced by a larger constant.\par 

Regarding the comparison of the volumes of $B_n$ and $G_n$, in Proposition 10 we will prove that 
$$|G_3|<|B_3|,\qquad |G_4|=|B_4|,\qquad |G_n|>|B_n|,\quad n\ge 5.\eqno(8)$$
The first two statements are easy to check, in fact $G_3$ is a proper subset of a translate of $B_3$ whereas $G_4$ is a translate of $B_4$. The last inequality in (8)  is not  so trivial to prove. Note also that $|G_n|\sim |B_n|$ as $n\to\infty$.
\medskip
With little or no extra effort we will  prove the following more general trace inequality:
\medskip\def\B{{\bar B}}
\proclaim Theorem 2. Let $\nu$ be a positive Borel measure on   $\B_n$ such that  for some $\lambda\in (0,n]$ and $C_0,r_0>0$
$$\nu\big(B(a,r)\cap \B_n)\le C_0 r^\lambda,\qquad \forall a\in\Rn,\;\;\forall r\in(0,r_0].\eqno(9)$$
There exists $C>0$ such that 
$$\int_{\B_n}\exp\bigg[{\lambda\alpha_n\over n}\,\bigg({|u(x)-Ru(x)|\over \|\Delta u\|_{n/2}^{}}\bigg)^{n\over n-2}\bigg]d\nu(x)\le C\eqno(10)$$
for all $u\in W^{2,n/2}(B_n),\;u\notin b^2$ where $\alpha_n$ is as in (6). The constant $\lambda\alpha_n/n$ in (10) is sharp provided there exists  $x_0\in \partial B_n$, if $n\ge4$, or  $x_0\in B_n$, if $n=3,4$, such that  $\nu(B(x_0,r)\cap\bar B_n)\ge C_1 r^\lambda$, for $0<r\le r_1$, some $C_1,r_1>0$. \par
\par\smallskip
Clearly Theorem 1 is a special case of Theorem 2, when $\nu$ is the Lebesgue measure. Another relevant special case is when $\nu=\H_{n-1}/S^{n-1}$, the $(n-1)-$dimensional Hausdorff measure restricted to the boundary of $B_n$. The result is the following boundary trace inequality:
\proclaim Corollary 3. There is $C>0$ such that 
$$\int_{S^{n-1}}\exp\bigg[{n-1\over n}\,\alpha_n\,\bigg({|u(x)-Ru(x)|\over \|\Delta u\|_{n/2}^{}}\bigg)^{n\over n-2}\bigg]d\H_{n-1}(x)\le C\eqno(11)$$
for all $u\in W^{2,n/2}(B_n),\;u\notin b^2$ where $\alpha_n$ is as in (6). The constant $(n-1)\alpha_n/n$ in (11) is sharp.\par

The proof of Theorem 2 is an application of results obtained by the authors in [FM]. We recall here the basic setup, in a simplified form.

Let  $(M,\mu)$,$\,(N,\nu)$ be  measure spaces with finite measure,  and suppose that $T$ is an integral operator of type
$$Tf(x)=\int_M K(x,z)f(z)d\mu(z)\,,\qquad x\in N\eqno(12)$$
where  $K:N\times M\to[-\infty,\infty]$ is measurable on  $(N\times M,\nu\times\mu)$. Define for $s>0$
$$\lambda_1(s,x)=\mu\big(\{z\in M: \,|K(x,z)|>s\}\big),\qquad x\in M$$
$$\lambda_2(s,z)=\nu\big(\{x\in N:\, |K(x,z)|>s\}\big),\qquad z\in N.$$
The following result was proved in [FM] (see [FM] Theorems 1 and 4   combined)
\proclaim Theorem 4 [FM]. In the above setup suppose that 
$$ \sup_{x\in M}\lambda_1(s,x)\le A s^{-\beta}\Big(1+O(\log^{-\gamma} s)\Big)\eqno(13)$$
and 
$$\sup_{z\in N}\lambda_2(s,z)\le B s^{-\beta_0}$$
as $s\to+\infty$,  for some $\beta,\gamma>1$, $\,0<\beta_0\le \beta$ and $A,B>0$. Then, $T$ is  defined by (12) on $L^{\beta'}(M)$ and 
there exists a constant $C$ such that 
$$\int_N \exp\bigg[{\beta_0\over  A\beta}\bigg({|Tf|\over \|f\|_{\beta'}}\bigg)^\beta\,\bigg]d\nu\le C\eqno(14)$$
for each $f\in L^{\beta'}(M)$, with $\ds{1\over \beta}+{1\over\beta'}=1.$ The constant $\beta_0/(A\beta)$ in (14) is sharp if the following additional conditions hold:\smallskip
\noin{i)} There is equality in (13)\smallskip
\noin{ii)} The supremum in (13) is attained at some $x_0\in M$\smallskip
\noin{iii)} There exists measurable sets $F_m\subseteq N,\, E_m\subseteq M$, $\,m\in \N$, such that for large $m$
$$E_m\supseteq\{y:\,|K(x_0,z)|>m\}\eqno(15)$$
$$\mu(E_m)\le C_1m^{-\beta}\eqno(16)$$
 $$C_2m^{-\beta_0}\le \nu(F_m)\le  C_3m^{-\beta_0}\eqno(17)$$ 
$$ \int_{M\setminus E_m} |K(x,z)- K(x_0,z)|\, |K(x_0,z)|^{\beta-1} d\mu(z)\le C_4\,,\qquad \forall x\in F_m
\eqno(18)$$
for some $C_1,C_2,C_3,C_4>0$. In particular, if 
$$\Phi_m(z)=K(x_0,z)|K(x_0,z)|^{\beta-2}\chi_{M\setminus E_m}^{}(z)$$ then $\Phi_m\in
 L^{\beta'}$, and 
$$\lim_{m\to\infty}\int_N \exp \bigg[\alpha\,\bigg({|T\Phi_m|\over\|\Phi_m\|_{\beta'}^{}}\bigg)^\beta\,\bigg]d\nu=+\infty,
\qquad \forall \alpha>{\beta_0\over A\beta}. $$
\smallskip
\smallskip

\par
The first step toward our proof of Theorem 2 is to write $u-Ru$ in terms of $\Delta u$, as an integral operator:
$$u(x)-Ru(x)=T(\Delta u)(x)=\int_{B_n}K(x,z)\Delta u(z)dz\qquad x\in \B_n$$
where $K(x,z)=N(x-z)-R(N(\cdot-z))$, where $N$ is the Newtonian potential as in (2). In other words, $T$ is the operator which gives what could be called ``the canonical solution" of the Poisson equation on the ball, that is the unique solution of $\Delta u=f$ which is orthogonal to the harmonic functions on the ball. It turns out that the kernel $K(x,z)$ can be explicitly computed, using well known formulas for the Bergman projection on the ball. In order state the precise result let us introduce some more terminology and notation.

\smallskip
For $x\in\Rn\setminus\{0\}$ let 
$$x^*={x\over |x|}$$
and define  the Dirichlet Green function as 
$$G(x,z)=N(x-z)-N(x^*-|x|z),\qquad z\in B_n,\;x\in\bar B_n,\;x\neq z$$
where $N(x-z)$ is the Newtonian potential as in (2) and (3), and with the convention that $N(0,z)=1$ if $z\neq 0$. It is well-known that $G$ is the fundamental solution of the Dirichlet problem on the ball.\smallskip
Define the extended Poisson Kernel (for $n\ge 2$) as 
$$P(x,y)={1\over\omega_{n-1}}{1-|x|^2|y|^2\over |x^*-|x|y|^n}\qquad x,y\in \R^n,\quad x^*\neq |x|y$$
so that if $y=y^*\in S^{n-1}$ then $P(x,y^*)$ is the standard Poisson kernel for the ball.

\bigskip
\proclaim Theorem 5. For each $p>1$ and each $f\in L^p(B_n)$, $n\ge2$,  the Poisson  equation $\Delta v=f$ has a unique solution $v\in (b^{p'})^\perp\cap W^{2,p}(B_n)$ given as $v=Tf$, where $T$ is the integral operator
$$Tf(x)=\int_{B_n} K(x,z)f(z)dz,\qquad x\in \bar B_n$$
with
$$K(x,z)=G(x,z)+{1-|z|^2\over2}\,P(x,z).\eqno(19)$$
Moreover, $T:L^p(B_n)\to (b^{p'})^\perp\cap W^{2,p}(B_n)$ is bounded and invertible, with inverse $\Delta$. 
\hfill\break
In particular, if $n\ge 3$ then $T:L^{n/2}(B_n)\to (b^{2})^\perp\cap W^{2,n/2}(B_n)$ is bounded and invertible, with inverse $\Delta$.
\par

The main novelty in the above theorem is the explicit formula for the kernel $K$, in (19). The formula looks simple enough for one to wonder whether it has appeared in print before; the autors were not able to find any published results of this sort in the literature. We should point out however an analogous, although unrelated, result in [HP] for the canonical solution of the  $\bar\partial$-Neumann problem on the unit ball of $\C^n$.
\medskip
The next step is to formulate and prove the following equivalent ``potential" version of Theorem 2:\medskip 
\proclaim Theorem 6. If $\nu$ is a positive Borel  measure on $\B_n$  as in (9), and  $T$ is the operator of Theorem 5, for $n\ge 3$, 
then there exists $C>0$ such that 
$$\int_{\B_n}\exp\bigg[{\lambda\alpha_n\over n}\,\bigg({|Tf(x)|\over \|f\|_{n/2}^{}}\bigg)^{n\over n-2}\bigg]d\nu(x)\le C\eqno(20)$$
for all $f\in L^{n/2}(B_n,dz)$. The constant $\lambda\alpha_n/n$ in (20) is sharp provided there exists  $x_0\in \partial B_n$, if $n\ge4$, or  $x_0\in B_n$, if $n=3,4$, such that  $\nu(B(x_0,r)\cap\bar B_n)\ge C_1 r^\lambda$, for $0<r\le r_1$, some $C_1,r_1>0$. \par
The proof of (20) follows at once from Theorem 4 and the following sharp distribution function  estimates for the kernel $K(x,z)$ defined in (19):

\proclaim Theorem 7.  
If 
$$\lambda_1(s,x)=|\{z\in B_n:\,|K(x,z)|>s\}|,\qquad x\in\B_n,$$
then there  exists $s_0>0$ and $\e>0$ such that for any $x\in\B_n$ and any $s\ge s_0$
$$\lambda_1(s,x)\le |B_n|\,c_n^{n\over n-2} s^{-{n\over n-2}}\big(1+O(s^{-\epsilon})\big),\qquad n=3,4\eqno(21)$$
$$\lambda_1(s,x)\le |G_n|\,c_n^{n\over n-2} s^{-{n\over n-2}}\big(1+O(s^{-\epsilon})\big),\qquad n\ge 5.\eqno(22)$$
For $n=3$ or $n=4$,  given any $x_0\in B_n$ (or any $x_0\in \p B_4$ if $n=4$)   we can choose $s_0$ so that equality occurs  in (21)  if $x=x_0$,  for any $s\ge s_0$. For $n\ge5$, given any $x_0\in \p B_n$ we can choose $s_0>0$ so that equality occurs in (22) if $x=x_0$, for any $s\ge s_0$.
\smallskip
If $n\ge 3$ and  $\nu$ is  a positive Borel  measure on $\B_n$  as in (9) and 
$$\lambda_2(s,z)=\nu\big(\{x\in \B_n:\,|K(x,z)|>s\}\big),\qquad  z\in B_n,$$
 then there exists $M>0$ such that for any $z\in B_n$
$$\lambda_2(s,z)\le M s^{-{\lambda\over n-2}},\qquad s> 0.\eqno(23)$$

\par

Estimate (23) will be easy to show, but the proofs of (21) and (22) -especially (22)- are surprisingly challenging. From Theorem 5 we know that the kernel $K$ is the sum of two kernels: the Dirichlet Green function $G(x,z)$   and the kernel ${1\over2}(1-|z|^2)P(x,z)$. Clearly $G(x,z)$ behaves like the Newtonian potential if $x$ is inside the unit ball, and it's 0 if $x$ is on the boundary. On the other hand, $P(x,z)$ is regular for $x$ inside the unit ball, but it becomes singular as $x$ approaches the boundary; in particular (see Lemma 5), for a boundary point $x^*\in \p B_n$ and $n
\ge 3$ 
$${1\over2}(1-|z|^2)P(x^*,z)=c_ng(x^*,z)\big|x^*-z\big|^{2-n}+O(|x^*-z|^{3-n}),\qquad x^*\in\partial B_n$$ where
$$g(x^*,z)=2(n-2)\,\bigg({x^*\cdot(x^*-z)\over |x^*-z|}\bigg)^2.$$
It is relatively easy to check that for {\it fixed $x$} inside the ball $\lambda_1(s,x)\sim|B_n|\,c_n^{n\over n-2} s^{-{n\over n-2}}$
as $s\to\infty$, and for $x$ on the boundary $\lambda_1(s,x)\sim|G_n|\,c_n^{n\over n-2} s^{-{n\over n-2}}$. The technical difficulty consists in establishing sharp asymptotic upper bounds for $\lambda_1(s,x)$  which are {\it uniform} with respect to $x\in \B_n$, by analyzing  carefully how the level sets and their measures change as $x$ moves towards the boundary,  due to the individual  contributions of the kernels $G$ and $P$ appearing in (19). Section 4 is dedicated to this analysis.

The sharpness statement of Theorem 6 will follow   from Theorem 4, the sharp estimates (21), (22), (23), and a smoothness H\"ormander-type estimate  on $K(x,z)$ (see (82) of Section~5). 
\medskip
A couple of remarks before concluding this section. In this paper we only treat the case of the unit ball, as our main domain. The main reason for that is that we have some known tools and explicit formulas for the Bergman projection at our disposal. It is natural to speculate that the results of this paper could be extended to any smooth domain. Even though on general smooth domains explicit formulas for the kernel $K$ as in Theorem 2 are in general hopeless, we speculate that near the boundary the behavior of $K$ should be similar to that of the ball kernel, so that one could try to adapt the arguments of this paper in the more general situation.

Finally, a few words should be spent regarding the case $n=2$. Obviously the exponential inequalities in Theorems 1,2,6 do not make any sense for $n=2$, nonetheless they can be replaced by different statements, in the same spirit as in a  result by  Brezis and Merle ([BM], Thm. 1), and the more recent results by Cassani, Ruf and Tarsi [CRT]. We will present these results in a forthcoming paper, as a special case of a more general class of exponential integral inequalities in the exceptional case where the dimension equals the order of the operator.

\vskip2em \centerline{\bf 3. Kernel computation: Proof of Theorem 5}\bigskip

We  recall the explicit formulas for the Bergman projection on the ball: for any $f\in L^2(B_n)$
$$Rf(x)=\int_{B_n}R(x,y)f(y)dy$$
where
$$R(x,y)={(n-4)|x|^4|y|^4+(8x\cdot y-2n-4)|x|^2|y|^2+n\over \omega_{n-1} |x^*-|x|y|^{n+2}}$$
is the reproducing kernel for the ball (see [ABR], Thm. 8.13, and [Li1], [Li2] but with a missing factor ${n-2\over2}$). 

We will actually find the following formula more useful:
$$R(x,y)=n P(x,y)+2{d\over dt}\bigg|_{t=1} P(tx,y)\eqno(24) $$
\bigskip
where 
$$P(x,y)={1\over\omega_{n-1}}{1-|x|^2|y|^2\over |y^*-|y|x|^n}={1\over\omega_{n-1}}{1-|x|^2|y|^2\over (1-2x\cdot y+|x|^2|y|^2)^{n/2}},\qquad x,y\in \R^n,\quad y^*\neq |y|x$$
is the extended Poisson kernel (so that if $y=y^*\in S^{n-1}$ then $P(x,y^*)$ is the standard Poisson kernel for the ball). This  formula is derived in [ABR], formula 8.12. We have that
$$P(ax,y)=P(x,ay),\qquad P(x,y)=P(y,x)=P(|x|y,x^*)$$
and that $P(x,\cdot)$ is harmonic on $\B_n$ for any $x\in B_n$.

 Assume for now that $f\in C^\infty(\B_n)$, and $n\ge 3$. The function 
$$g(x)=N*f(x)=-\int_{B_n}c_n|x-z|^{2-n}f(z)dz$$
is $C^\infty(\B_n)$ and solves $\Delta g=f$ on $B_n$  (and $\Delta g=0$ on $\B_n^c$), so $g-Rg$ is the unique function $v$ of $L^2(B_n)$ such that $\Delta v=f$ and $v\in (b_2)^\perp$.  Therefore, it's enough to compute
$$H(x,z):=c_n\int_{B_n}R(x,y)|y-z|^{2-n}dy,\qquad x\in \B_n.$$
Using formula (24)
$$H(x,z)=n c_n \int_B P(x,y)|y-z|^{2-n}dy+2c_n{d\over dt}\bigg|_{t=1} \int_B P(tx,y)|y-z|^{2-n}dy\eqno(25)$$
First write
$$\eqalign{\int_B P(x,y)|y-z|^{2-n}dy&=\int_0^1 r^{n-1}dr\int_{S^{n-1}}P(x,ry^*)|ry^*-z|^{2-n}dy^*=\cr&=\int_0^1 rdr \int_{S^{n-1}}P(rx,y^*)|y^*-zr^{-1}|^{2-n}dy^*.\cr}\eqno(26)$$

Now for given $0<r<1$
$$\int_{S^{n-1}}P(x,y^*)|y^*-zr^{-1}|^{2-n}dy^*=\cases{|x-zr^{-1}|^{2-n} & if $|z|>r$\cr \cr \big|x^*-|x|zr^{-1}\big|^{2-n} & if $|z|<r$\cr}$$
since the functions on the right are harmonic and with the same boundary values as the function on the left. Hence, evaluating the above formulas  at $rx$ and inserting them in (26) yields
$$\eqalign{\int_B P(x,y)|y-z|^{2-n}dy&=\int_0^{|z|}r|rx-zr^{-1}\big|^{2-n}dr+\int_{|z|}^1 r \big|x^*-|x|z\big|^{2-n}dr=\cr
&=\int_0^{|z|}r|rx-zr^{-1}\big|^{2-n}dr+{1\over2}(1-|z|^2)\big|x^*-|x|z\big|^{2-n}\cr}$$
and
$$\eqalign{&{d\over dt}\bigg|_{t=1} \int_B P(tx,y)|y-z|^{2-n}dy\cr&={d\over dt}\bigg|_{t=1}
\int_0^{|z|}r|rtx-zr^{-1}|^{2-n}dr+{1\over2}(1-|z|^2){d\over dt}\bigg|_{t=1}
\big|x^*-t|x|z\big|^{2-n}\cr&={d\over dt}\bigg|_{t=1}
\int_0^{|z|\sqrt t}t^{-n/2} r|rx-zr^{-1}|^{2-n}dr+{n-2\over2}(1-|z|^2)|x|\,{z\cdot(x^*-|x|z)\over\big|x^*-|x|z\big|^n}\cr
&={|z|^2\over2}\,\big|x|z|-z^*\big|^{2-n}-{n\over2}\int_0^{|z|} r|rx-zr^{-1}|^{2-n}dr+{n-2\over2}(1-|z|^2)|x|\,{z\cdot(x^*-|x|z)\over\big|x^*-|x|z\big|^n}\cr}$$

Putting all this in  (25) and using $|x^*-|x|z|=|z^*-|z|x|$
gives the explicit formula
$$\eqalign{H(x,z)&=c_n|z|^2\,\big|x^*-|x|z\big|^{2-n}+c_n(n-2)(1-|z|^2)|x|\,{z\cdot(x^*-|x|z)\over\big|x^*-|x|z\big|^n}\cr&+{nc_n\over2} (1-|z|^2)\big|x^*-|x|z\big|^{2-n},\cr}$$
and after a few more simple algebraic calculations (19) is obtained. The the case $n=2$ (and $f\in C^\infty(\B_n)$) is derived similarly, using $-(2\pi)^{-1}\log|x-z|$ in place of $c_n|x-z|^{2-n}$, and a few minor changes in the proof.
\smallskip
The proof of Theorem 5 is completed by observing that the operator $f\to N\*f$ is bounded from $L^p$ to $W^{2,p}$ of the ball, for $p>1$ (see for ex. [GT], Thm 9.9), and the operator $R$ is bounded from $W^{2,p}$ to itself (see for example [KK], Thm. 4.5). Hence the operator $f\to Tf=N*f-R(N*f)$ is bounded from $L^p$ to $W^{2,p}$. The fact that $Tf\in (b^{p'})^\perp$ follows by a density argument.\bigskip

\eject

\centerline{\bf 4. Kernel distribution estimates: proof of Theorem 7}\bigskip

 For simplicity we will work with the normalized kernel $$K_0(x,z)=-c_n^{-1}K(x,z)=|x-z|^{2-n}-|x^*-|x|z|^{2-n}-{n-2\over2}\,{(1-|z|^2)(1-|x|^2|z|^2)\over |x^*-|x|z|^n}.$$
From now on, and with a slight abuse of notation, we will let 
$$\lambda_1(s,x)=|\{z\in B_n:\,|K_0(x,z)|>s\}|,\qquad s>0,\;x\in\B_n$$
$$\lambda_2(s,z)=\nu\big(\{x\in\B_n:\,|K_0(x,z)|>s\}\big),\qquad s>0,\;z\in B_n$$
 where $\nu$ is a Borel measure on $\B_n$ satisfying (9). 
Note that $\lambda_1$ is invariant under rotations.

The inequalities in Theorem 7 are equivalent to the following: 
$$\lambda_1(s,x)\le |B_n|\, s^{-{n\over n-2}}\big(1+O(s^{-\epsilon})\big),\qquad n=3,4\eqno(27)$$
$$\lambda_1(s,x)\le |G_n|\, s^{-{n\over n-2}}\big(1+O(s^{-\epsilon})\big),\qquad n\ge 5,\eqno(28)$$
valid for $s\ge s_0$ uniformly in $x\in \B_n$, and 
$$\lambda_2(s,z)\le M s^{-{\lambda\over n-2}},\qquad s> 0.\eqno(29)$$
uniformly in $z\in B_n$.

\bigskip
The proofs of the above inequalities are divided in six main  steps:

\smallskip
\noin Step 1: we derive  an asymptotic expansion of $K_0$  around its singularities
\smallskip\noin
Step 2:  we easily prove (29). 
\smallskip\noin
Step 3: we prove (27) and (28) with the equality sign when $x$ is on the boundary
\smallskip\noin
Step 4: we prove that $|B_3|>|G_3|,\,|B_4|=|G_4|,\,$ and $|B_n|< |G_n|$ for $n\ge 5$.
\smallskip\noin
Step 5: we prove the inequality in (27) for all $n$, uniformly  in the range $|x|\le 1-s^{-{1\over n-2}}$
\smallskip\noin
Step 6: we prove (27) and (28) uniformly  in the range $ 1-s^{-{1\over n-2}}\le |x|\le 1.$

\medskip\noin{\bf {Convention: }} {\it Throughout the paper $\e$, $C$, $s_0$, and $t_0$ will denote suitable positive constants depending at most on the dimension $n$. Such constants  might take different values  even within a single chain of identities or inequalities, and their  precise values is irrelevant for our purposes.  }
\medskip
\bigskip\eject
\noindent{\it\underbar{Step 1: Kernel Asymptotics}}\bigskip
Although the kernel $K_0$ is a difference of two good looking  positive kernels, in order to compute the asymptotics of $\lambda_1$ and $\lambda_2$ we find it  more useful to just deal with the following asymptotic and global  estimates of $K_0$:

\proclaim Lemma 8. The following asymptotic expansions hold for any $b\in (0,1)$:\smallskip 
\item{a)} If $|x|\le b<1,\;z\in B_n$  and $x\neq z$  
 then 
$$K_0(x,z)=|x-z|^{2-n}+O(1)\eqno(30)$$
\smallskip\item {b)} If $b\le |x|\leq 1,\;z\in B_n$ and $x\neq z$  then 
$$K_0(x,z)=|x-z|^{2-n}-\big(1+g(x,z)\big)\big|x^*-|x|z\big|^{2-n}+O(|x-z|^{3-n})\eqno(31)$$ where
$$g(x,z)=2(n-2)\,{x^*\cdot(x^*-z)\,x^*\cdot (x^*-|x|z)\over |x^*-|x|z|^2}\ge0$$
and where the $O'$s are  uniform in the respective domains of $(x,z)$. 
\smallskip\noin
The following global estimates hold  for $ x\in \bar B_n,\;z\in B_n,\;x\neq z$:
$$K_0(x,z)\ge -g(x,z)\big|x^*-|x|z\big|^{2-n}\eqno(32)$$
and
$$|K_0(x,z)|\le H |x-z|^{2-n},\qquad\eqno(33)$$
for some $H>0$ independent of $x,z$.\par
 \bigskip\pf Proof.
The asymptotic expansion in (31) follows easily from the following facts:\smallskip \smallskip\item {1.} $H(x,z)$ is continuous  on $\overline{B(0,b)}\times \overline{B(0,1)}$ for any $b<1$; \smallskip\item{ 2.}
  $1-|x|^2|z|^2=2x^*\cdot(x^*-|x|z)-\big|x^*-|x|z\big|^2$;\smallskip
\item{3.} $\big||x|-1\big|= |x-x^*|\le |x-z|+|z-x^*|$; \smallskip
\item{4.} $|x-z|\le |x^*-|x|z|$;\smallskip\item{5.} $|x^*-z|\le b^{-1} |x^*-|x| z|$ if $b\le |x|\le 1$ and $|z|\le 1$.
\smallskip\noin
Finally, (32) follows from fact {2.} \hskip-.05  em above, and (33) is a simple consequence of (30) and (31).\endpf
\bigskip
\bigskip\eject\noindent{\it\underbar{Step 2: Global estimate on $\lambda_2$}}\bigskip
\medskip
Inequality (29) follows immediately from (33) and the assumptions on $\nu$: for each $z\in B_n$ and $s>0$
$$\lambda_2(s,z)=\nu\big(\{x\in\B_n:\,|K_0(x,z)|>s\}\big)\le
\nu\big(\{x\in\B_n:\,|x-z|<(s/H)^{-{1\over n-2}}\}\big)\le C s^{-{\lambda\over n-2}}.$$
\bigskip\noindent{\it\underbar{Step 3: Asymptotics  of $\lambda_1$ on the boundary}}\bigskip
\proclaim Proposition 9. There exists $s_0>0$ such that for any $x^*\in \p B_n$ we have
$$\lambda_1(x^*,s)=|G_n| s^{-{n\over n-2}}\Big(1+O(s^{-\e})\Big),\qquad s\ge s_0\eqno(34)$$
and 
$$|G_n|={\pi^{n-1\over2}\over n}(2n-4)^{n\over n-2}{\Gamma\big({1\over2}+{n\over n-2}\big)\over\Gamma\big({n\over2}+{n\over n-2}\big)}.\eqno(35)$$\par
\medskip
\pf Proof. We have
$$K_0(x^*,z)=-(2n-4)\,\bigg({x^*\cdot(x^*-z)\over|x^*-z|}\bigg)^2|x^*-z|^{2-n}+O\big(|x^*-z|^{3-n}\big).$$
By rotation invariance we can assume 
$$x^*=e_1=(1,0,0,...,0).$$
 and let's also let $$w=x^*-z,\qquad g(w^*)=(2n-4)(w_1^*)^2, \qquad w\in\R^n.$$
With this notation 
$$G_n=\{w:\,|w|^{n-2}\le g(w^*)\}.$$
and
$|K(e_1,z)|\le C|w|^{3-n}$, so that if $|K(x^*,s)|>s$ then $|w|\le C s^{-\e}$ (with $\e=-\un$). Hence
$$\eqalign{\lambda_1(e_1,s)&\le |\{w: |x^*-w|<1,\, |w|^{2-n}\big(g(w^*)+C|w|\big)>s\}|\cr& \le|\{w: \, |w|^{2-n}\big(g(w^*)+C s^{-\e}\big)>s\}|. \cr}$$
After passing in polar coordinates  it's easy to check that
$$\lambda_1(e_1,s)\le{1\over n}\int_{S^{n-1}}
\bigg({g(w^*)\over s}\bigg)^{n\over n-2}dw^*+Cs^{-{n\over n-2}-\e},$$
and 
$$|G_n|={1\over n}\int_{S^{n-1}}
\big(g(w^*)\big)^{n\over n-2}dw^*.$$
Using the formula
$$\int_{S^{n-1}} F(w_1^*)dw^*=\omega_{n-2}\int_{-1}^1 F(t)(1-t^2)^{n-3\over2}dt$$
we get 
$$\eqalign{\int_{S^{n-1}}
\big(g(w^*)\big)^{n\over n-2}dw^*&={\omega_{n-2}\over n}(2n-4)^{n\over n-2}\int_0^1 t^{2n\over n-2}(1-t^2)^{n-3\over2}dt=\cr&={\omega_{n-2}\over 2n}(2n-4)^{n\over n-2}\int_0^1 t^{n+2\over2n-4}(1-t)^{n-3\over2}dt=\cr
&={\omega_{n-2}\over 2n}(2n-4)^{n\over n-2}{\Gamma\big({3n-2\over2n-4}\big)\Gamma\big({n-1\over2}\big)\over\Gamma\big({n^2\over 2n-4}\big)}\cr
}$$ 
which proves (35) and also  (34), but with ``$\le$". To derive  equality in (34) we note that
$$|K(e_1,z)|\ge |w|^{2-n}\big(g(w^*)-C|w|\big) $$
and that 
$$|w|^{2-n}\bigg[g(w^*)-C\bigg({g(w^*)\over s}\bigg)^\un\bigg]>s\imp |w|^{2-n}\big(g(w^*)-C|w|\big)>s.$$
Thus
$$\eqalign{\lambda_1(e_1,s)&\ge\bigg|\bigg\{w:\,|e_1-w|<1,\,|w|^{2-n}\bigg[g(w^*)-C\bigg({g(w^*)\over s}\bigg)^\un\bigg]>s\bigg\}\bigg|\cr&\ge \Big|\Big\{w:\,|w|^2<2w_1^*,\,|w|^{2-n}\Big(g(w^*)-Cs^{-\un}\Big)>s\Big\}\Big|\cr
&=\bigg|\bigg\{w: g(w^*)\ge C s^{-\un},\,|w|\le\min\Big\{\sqrt {2w_1^*},\Big(g(w^*)-Cs^{-\un}\Big)^{\un} s^{-\un}\Big\}\bigg\}\bigg|\cr
}$$
and it's easy to check that if $g(w^*)>Cs^{-\un}$ then 
$\sqrt {2w_1^*}>\Big(g(w^*)-Cs^{-\un}\Big)^{\un} s^{-\un}$ for $s$ large enough. Hence
$$\eqalign{\lambda_1(e_1,s)&\ge\bigg|\bigg\{w: g(w^*)\ge C s^{-\un},\,|w|\le\Big(g(w^*)-Cs^{-\un}\Big)^{\un} s^{-\un}\bigg\}\bigg|\cr
}$$
and if $H_s=\Big\{w\in S^{n-1}:\,g(w^*)>Cs^{-\un}\Big\}$
then 
$$\eqalign{\lambda_1(e_1,s)&\ge{s^{-{n\over n-2}}\over n}\int_{H_s}\Big(g(w^*)-Cs^{-\un}\Big)^{n\over n-2}dw^*
\ge{s^{-{n\over n-2}}\over n}\int_{H_s}\Big[\big(g(w^*)\big)^{n\over n-2}-Cs^{-\e}\Big]dw^*\cr&\ge {s^{-{n\over n-2}}\over n}\int_{S^{n-1}}\big(g(w^*)\big)^{n\over n-2}dw^*-Cs^{{-{n\over n-2}}-\e}\cr}$$\endpf
\bigskip\eject
\noindent{\it\underbar{Step 4: Comparing the volumes of $B_n$ and $G_n$} }\bigskip
We now give a comparison theorem for the volumes of $B_n$ and $G_n$.
The result does not appear to be provable using trivial or straightforward  methods, such as induction.\medskip
\proclaim Proposition 10. The following hold:
\medskip
\item{a)} $\;|B_3|>|G_3|$\smallskip
\item{b)} $\;|B_4|=|G_4|$\smallskip
\item{c)} $\;|B_n|<|G_n|$, for $n\ge 5$.
\par
\medskip\pf Proof. Recall
$$|B_n|={2\pi^{n/2}\over n \Gamma\big({n\over2}\big)},\qquad|G_n|={\pi^{n-1\over2}\over n}(2n-4)^{n\over n-2}{\Gamma\big({1\over2}+{n\over n-2}\big)\over \Gamma\big({n\over2}+{n\over n-2}\big)}$$

We then have

$$|B_3|={4\pi\over3}>|G_3|={16\pi\over 21}$$
$$|B_4|=|G_4|={\pi^2\over 2}.$$
The inequality in c) is equivalent to 
 $${\pi^{1/2}\over \Gamma(n/2)}< {1\over2} (2n-4)^{n\over n-2}\,{\Gamma\big({n\over n-2}+{1\over2}\big)\over\Gamma\big({n\over n-2}+{n\over 2}\big)},\qquad n\ge5\eqno(36)$$

Letting  $t={2\over n-2}\in(0,1]$ inequality  (36) becomes

$${\Gamma\big(1+{1\over t}\big)\over\Gamma\big(t+{1\over t}\big)}\,{\Gamma\big({1\over2}+t+1\big)\over\Gamma\big(1+{1\over2}\big)} {2^{2t} t^{-1-t}\over \big(t+1+{1\over t}\big)\big(t+{1\over t}\big)}\ge 1$$

Using the inequality  (see [Ke])
$${\Gamma(x+1)\over\Gamma(x+\lambda)}>
\bigg(x+{\lambda\over2}\bigg)^{1-\lambda},\qquad x>0,\,\lambda\in(0,1)\eqno(37)$$
we obtain 
$${\Gamma\big({1\over t}+1\big)\over\Gamma\big({1\over t}+t\big)}>\bigg({1\over t}+{t\over 2}\bigg)^{1-t}$$
and using again (37) but with $1-\lambda=t$ and $x=\half+t$

$${\Gamma\big({1\over2}+t+1\big)\over\Gamma\big(1+{1\over2}\big)}  > \bigg(1+{t\over 2}\bigg)^t$$
so that the left hand side of (36) is greater than 
$$K(t)=\bigg({1\over t}+{t\over 2}\bigg)^{1-t} \bigg(1+{t\over 2}\bigg)^t  {2^{2t} t^{-1-t}\over \big(t+1+{1\over t}\big)\big(t+{1\over t}\big)}=2^{-1+2t} \,{(t^2+2)^{1-t}(t+2)^t\over (t^2+t+1)(t^2+1)}.$$

Note that  $K(0)=K(1)=1$, so it is enough to show that if $H(t)=\log K(t)$ then $H''(t)< 0$, on $[0,1]$; this can be checked by a straightforward (but lenghty) algebraic calculation, which shows that $-H''(t)$ is a ratio of two polynomials with positive coefficients.\endpf

\bigskip\noindent{\it\underbar{Step 5: uniform estimates on $\lambda_1$ in the range $|x|\le 1-s^{-{1\over n-2}}$} }\bigskip
\medskip
\proclaim Proposition 11.   For any $n\ge 3$ there is $s_0>0$ and $\e>0$ so that for $|x|\le 1- s^{-{1\over n-2}}$
$$\lambda_1(s,x)\le  |B_n|\,  s^{-{n\over n-2}}\big(1+O\big(s^{-\e}\big)\big),\qquad s\ge s_0. \eqno(38)$$
Given any $x_0\in B_n$ we can choose $s_0$ such that equality occurs in (38) when $x=x_0$.

\par
\pf Proof.  
Let us first show (38) in the easier case $|x|\le b<1$, for any given $b$ with $0< b<1$. This follows from (30):
$$\lambda_1(s,x)\le|\{z\in B_n: |x-z|^{2-n}+C>s\}|\le |B_n|(s-C)^{-{n\over n-2}}=|B_n| s^{-{n\over n-2}}(1+O(s^{-\e})),$$
for $s\ge s_0$. If $x_0\in B_n$ is given, choose $b$ so that $|x_0|<b<1$ and  one can reverse the above inequality when $x=x_0$  in a similar way:
$$\lambda_1(s,x_0)\ge|\{z\in B_n:|x-z|<(C+s)^{-{1\over n-2}}\}|=|B_n|(C+s)^{-{n\over n-2}}=|B_n| s^{-{n\over n-2}}(1+O(s^{-\e})),$$
provided $(C+s)^{-{1\over n-2}}<b-|x_0|$.
\medskip
Suppose now that $0<b<1$ and  $b\le |x|\le 1-s^{-\un}$, for $s\ge s_0$ large enough, and let us analyze in  more detail
the sets $\{z: K_0(x,z)>s\}$ and $\{z: -K_0(x,z)>s\}$ under these assumptions.
First, note that 
$$K_0(x,z)\le |x-z|^{2-n}+C|x-z|^{3-n}\le H|x-z|^{2-n}$$
so 
$$\eqalign{|\{z:K_0(x,z)>s\}|&\le \{z: |x-z|^{2-n}+C s^{-{1\over n-2}}|x-z|^{2-n}>s\}|\cr& \le |B_n| s^{-{n\over n-2}}\Big (1+Cs^{-{1\over n-2}}\Big)^{n\over n-2}=|B_n| s^{-{n\over n-2}}\big(1+O(s^{-\e})\big).\cr}$$  

Next, we have 
$K_0(x,z)\ge - g(x,z)|x^*-|x|z|^{2-n}$ (see (32))
so that
$$|\{z: -K_0(x,z)>s\}|\le |\{z:g(x,z)|x^*-|x|z|^{2-n}>s\}|\eqno(39)$$
and we now claim that the right-hand side of (39)  is actually 0 for $|x|\le 1-s^{-{1\over n-2}}$.
Note that $|x^*-|x|z|\ge 1-|x|(x^*\cdot z)$, so
$$g(x,z)|x^*-|x|z|^{2-n}\le (2n-4)|x|^{1-n}(1-x^*\cdot z)\big(|x|^{-1}-x^*\cdot z\big)^{1-n}.$$
The function 
$$\varphi(\zeta)=(1-\zeta)(|x|^{-1}-\zeta)^{1-n},\qquad \zeta \in [-1,1]$$
attains a global maximum at 
$$\zeta=1-{|x|^{-1}-1\over n-2}\in [-1,1]$$
for ${1\over 2n-3}\le |x|\le 1$
 and so
$$\eqalign{{(2n-4)\over|x|^{n-1}}(1-\zeta)\big(|x|^{-1}-\zeta\big)^{1-n}&\le {2(1-|x|)\over |x|^n}\bigg(|x|^{-1}-1+{|x|^{-1}-1\over n-2}\bigg)^{1-n}\cr&= {2\over|x|}(1-|x|)^{2-n}\bigg(1-{1\over n-1}\bigg)^{n-1}\le {2\over e |x|}\,(1-|x|)^{2-n}.\cr}$$
Hence the set $\{z\in B_n: g(x,z)|x^*-|x|z|^{2-n}>s\}$ is empty if $2|x|^{-1}(1-|x|)^{2-n}e^{-1}\le s$, and this is certainly true if $2/e\le |x|\le 1-s^{-{1\over n-2}}$, and  in particular if $b$ is chosen so that $2/e\le b<1$. This settles (38). To conclude the proof, if $x_0\in B_n$ we can choose $b$ so that $\max\{|x_0|,2/e\}<b<1$ and the previous discussion guarantees that (38) can be reversed when $x=x_0$, for $s\ge s_0$ large enough.

\endpf\vskip2.7em 
\noindent{\it\underbar{Step 6: uniform estimates on $\lambda_1$ in the range $1-s^{-{1\over n-2}}\le |x|\le 1$}}\bigskip
The next task, and the most challenging one,  is to analyze $\lambda_1(s,x)$ in the range $1-s^{-{1\over n-2}}\le|x|\le 1$.  In particular we want to prove that for some $s_0>0$ the following estimates hold for $s\ge s_0$ and $1-s^{-{1\over n-2}}\le|x|\le 1$
$$\lambda_1(s,x)\le |B_n|\,s^{-{n\over n-2}}\big(1+O(s^{-\e})\big),\qquad n=3,4\eqno(40)$$
$$\lambda_1(s,x)\le |G_n|\,s^{-{n\over n-2}}\big(1+O(s^{-\e})\big),\qquad n\ge5\eqno(41)$$
which settle completely (27) and (28) and therefore Theorem 7.\medskip\eject
We begin by observing that the condition $|K_0(x,z)|>s$ together with (33)  implies $|x-z|\le Cs^{-\un}$,   and, as a consequence,   the condition  $1-s^{-{1\over n-2}}\le|x|\le 1$ together with  (31)  implies
$$K_0(x,z)\le |x-z|^{2-n}-(1+g(x,z))|x^*-|x|z|^{2-n}+ C s^{-\un}|x-z|^{2-n}.\eqno(42)$$

By rotation invariance we can assume that 

$$x=x_1e_1,\qquad e_1=(1,0,0..,0)\in \R^n,$$

and we  make the following  convenient change of variables:
$$t=s^{-\un}\le t_0<1\qquad x=(1-\theta t)e_1, \,\,\qquad z=e_1-t y,\qquad \,0\le \theta\le 1.\eqno(43) $$
We also let 
$$y=(y_1,y_2,...,y_n)=(y_1,y')\in\Rn.$$
With  this notation, and using (42) and (32),  we get
$$\{z\in B:\,K_0(x,z)>s\}\subseteq e_1-s^{-\un} E(\theta,t)$$
$$\{z\in B:\,-K_0(x,z)>s\}\subseteq e_1-s^{-\un} D_0(\theta)\subseteq e_1-s^{-\un} D(\theta)$$
where
$$E(\theta,t):=\Big\{y: y_1\ge0,\,(1+C t)|y-\theta e_1|^{2-n}-\bigg(1+(2n-4){y_1(y_1+\theta-\theta t y_1)\over|y+\theta e_1-\theta t y|^2}\bigg)\,|y+\theta e_1-\theta t y|^{2-n}>1\Big\}$$
$$D_0(\theta)=\Big\{y: y_1\ge0 ,\,\bigg(1+(2n-4){y_1(y_1+\theta)\over|y+\theta e_1|^2}\bigg)\,|y+\theta e_1|^{2-n}-|y-\theta e_1|^{2-n}>1\Big\}\eqno(44)$$
$$D(\theta)=\{y\in\R^n:\, |y+\theta e_1|^n\le (2n-4)y_1(y_1+\theta)\}\supseteq D_0(\theta).$$
Let us also define
$$E(\theta)=E(\theta,0)=\Big\{y: y_1\ge0 ,\,|y-\theta e_1|^{2-n}-|y+\theta e_1|^{2-n}\bigg(1+(2n-4){y_1(y_1+\theta)\over|y+\theta e_1|^2}\bigg)>1\Big\}.$$

\medskip 

\proclaim Proposition 12. There exist $s_0,\,\e>0$ such that for $s\ge s_0$
$$\lambda_1(s,x)\le s^{-{n\over n-2}}|E(\theta)\cup D(\theta)|\,\big(1+O(s^{-\e})\big),$$
for all  $x$ and $\theta$  related as in (43).
\par\eject
\pf Proof of Proposition 12. Since 
$$\lambda_1(s,x)\le s^{-{n\over n-2}}|E(\theta,t)\cup D(\theta)|$$
the proof is completed once we show that for some  $C,\e,t_0>0$ 
$$|E(\theta,t)\setminus E(\theta)|\le Ct^\e,\qquad t\le t_0,\quad 0\le \theta\le 1.\eqno(45)$$

 We begin with the following inequalities 
$$|y+\theta e_1|^{2-n}\le |y+\theta-t\theta y|^{2-n}$$
$${y_1(y_1+\theta)\over|y+\theta e_1|^2}\le {y_1(y_1+\theta-t\theta y_1)\over|y+\theta e_1-t\theta y|^2},$$
valid for  $y_1\ge0$ and $t\le t_0$, and whose proof is straightforward.

As a consequence we get the following: let 
$$F_\theta(y)=|y-\theta e_1|^{2-n}-|y+\theta e_1|^{2-n}\bigg(1+(2n-4){y_1(y_1+\theta)\over|y+\theta e_1|^2}\bigg)$$
then 
$$E(\theta,t)\subseteq \{y:\,y_1\ge0,\, F_\theta(y)\ge 1-Ct|y-\theta e_1|^{2-n} \}\eqno(46)$$

\medskip
and clearly 
$$E(\theta)=E(\theta,0)=\{y:\,y_1\ge0,\, F_\theta\ge 1\}.$$

If $|y-\theta e_1|\ge t^{1\over 2n-4}$, then $t|y-\theta e_1|^{2-n}\le \sqrt t$,  hence
$$E(\theta,t)\setminus E(\theta)\subseteq  \{y: y_1\ge0, \, F_\theta(y)\ge 1-C\sqrt t\,\}\cup \{y:|y-\theta e_1|\le t^{1\over 2n-4}\},$$
 and since $F$ is invariant under rotations about the $y_1-$ axis, it suffices to prove 
$$|\{(y_1,y_2):\,y_1,y_2>0,\, 1-C\sqrt t\le F_\theta(y_1,y_2)\le 1\}|\le C t^\e,\quad t\le t_0,\quad 0\le \theta\le 1.\eqno(47)$$
where (with a slight  abuse of notation) $F_\theta(y_1,y_2)$ is the section of $F_\theta(y)$ on the plane \break $y_3=...=y_n=0$. Note also that $F_\theta(y_1,y_2)$ is well defined and smooth in the region $y_1,y_2>0$, for any    
 $\theta\in[0,1]$.

\proclaim Lemma 13. Given any $a>0$, if $F_\theta(y_1,y_2)\ge a$ for some $y_1,y_2>0$ then, with $y=(y_1,y_2)$,
$${\p F_\theta\over \p y_2}(y)\le -{any_2\over|y+\theta e_1|^2}.\eqno(48)$$\par\smallskip
\pf Proof of Lemma 13. 
We have 
$$\eqalign{{\p F_\theta\over \p y_2}&=y_2(2-n)|y-\theta e_1|^{-n}
-y_2(2-n)|y+\theta e_1|^{-n}\bigg(1+(2n-4){y_1(y_1+\theta)\over|y+\theta e_1|^2}\bigg)+\cr&\hskip5em+
|y+\theta e_1|^{2-n}(2n-4){2y_2 y_1(y_1+\theta)\over|y+\theta e_1|^4}\cr&=y_2(n-2)|y+\theta e_1|^{-n}
\bigg[-{|y+\theta e_1|^{n}\over|y-\theta e_1|^{n}}+1+(2n-4){y_1(y_1+\theta)\over|y+\theta e_1|^2}+4{y_1(y_1+\theta)\over|y+\theta e_1|^2}\bigg]\cr&=
y_2(n-2)|y+\theta e_1|^{-n}
\bigg[-\bigg({|y+\theta e_1|\over|y-\theta e_1|}\bigg)^n+1+2n{y_1(y_1+\theta)\over|y+\theta e_1|^2}\bigg].\cr}
$$
If  $F_\theta(y_1,y_2)\ge a$ we obtain 
$$\bigg({|y+\theta e_1|\over|y-\theta e_1|}\bigg)^{n-2}-1-(2n-4){y_1(y_1+\theta)\over|y+\theta e_1|^2}\ge a\, |y+\theta e_1|^{n-2}$$
and letting 
$$R={|y+\theta e_1|\over|y-\theta e_1|}\ge1$$
gives that 
$$\eqalign{{\p F_\theta\over \p y_2}&\le 
y_2(n-2)|y+\theta e_1|^{-n}
\bigg[-R^n+1+{n\over n-2}\Big(R^{n-2}-1-a|y+\theta e_1|^{n-2}\Big)\bigg]\cr&
=y_2|y+\theta e_1|^{-n}
\Big(-(n-2)R^n-2+n R^{n-2}-an|y+\theta e_1|^{n-2}\Big).\cr}$$
The function $(n-2)R^n-n R^{n-2}+2$ has a minimum at $R=1$, where it vanishes, and (48) is proved.
\endpf

Lemma 13 easily implies that for each $a>0$ and each fixed $y_1>0$ the vertical section $\{y_2>0: F_\theta(y_1,y_2)\ge a\}$ is either the empty set or a vertical segment $\{(y_1,v),\, v\in (0,w]\,\}$, some $w=w(\theta,y_1,a)>0$. Indeed, if $F(y_1,v)\ge a$ for some $v>0$ and $F(y_1,u)<a$ for some $u\in (0,v)$, then by continuity of $F_\theta$ we can find a smallest $v^*>u$ such that $F_\theta(y_1,v^*)\ge a$. But Lemma 13 guarantees that $\partial_{y_2}F_\theta(y_1,v^*)<0$, and this contradicts the minimality of~$v^*$, since  $F_\theta(y_1,v)>F(y_1,v^*)$,  for some $v\in (u,v^*).$ 

Also, if $F_\theta(y)\ge a$ then $|y-\theta e_1|^{2-n}\ge a$ and so $|y-\theta e_1|\le a^{-\un}$,  which means that the level set $\{y:\,F_\theta(y)\ge a\,\}$ is inside the ball of radius $1+a^{-\un}$, and all of its nonempty vertical sections in the first open quadrant must be bounded, half-closed segments.

Taking $a=1-C\sqrt t\ge {1\over 2}$ for $t\le t_0$, we obtain that the set $\{(y_1,y_2):\,y_1,y_2>0,\, 1-C\sqrt t\le F_\theta(y_1,y_2)\le 1\}$ is inside a ball or radius 3, and its vertical sections are either contained in the strip 
$\{ 0\le y_1\le 3,\, 0\le y_2\le  t^{1/4}\}$, or else they are segments with length smaller than 
$${1-(1-C\sqrt t)\over\ds{n t^{1/4}\over 32}}\le C t^{1/4}$$
since along those segments (48) implies 
$$ {\p F_\theta\over\p y_2}(y_1,y_2)\le-{{1\over2}nt^{1/4}\over 16}.$$
 From these results (47) follows easily, and hence Proposition 12 is  proved.\endpf
\smallskip
\proclaim Proposition 14. For $\theta\in[0,1]$ we have
$$|E(\theta)\cup D(\theta)|\le\cases{|B_n| & if $n=3,4$\cr  \cr|G_n| & if $n\ge5.$\cr}\eqno(49)$$
\par

Once this is done, inequalities (40) and (41) and Theorem 7 are completely proved.\smallskip

\pf Proof of Proposition 14. The strategy of this proof  is to first show that  for certain ranges of $\theta$ the sets $E(\theta)$ and $D(\theta)$ are either both inside the ball $\bar{B}_n+\theta e_1$ or both inside  $G_n-\theta e_1(\theta)$. Unfortunately, however, it does not seem possible to argue with inclusions for 
all values of $\theta$ in the interval $[0,1]$, and for all values of $n$; a critical range of $\theta'$s exists for which  the inequality in (49) for $n\ge 6$ will be  proved by actually estimating certain integrals.
\medskip
Define
$$B(\theta)=\overline{B}_n+\theta e_1=\{y:\,|y-\theta|\le1\}$$
$$G(\theta)=\overline G_n-\theta e_1=\{y:\,0\le y_1+\theta\le (2n-4)^{\un},\; |y+\theta e_1|^n\le (2n-4)(y_1+\theta)^2\},$$
which are obtained by rotating the regions under the curves $y_2=h^{1/2}(y_1-\theta)$ and \break $y_2=f^{1/2}(y_1+\theta)$ where
$$h(v)=1-v^2,\qquad f(v)=(2n-4)^{2/n}v^{4/n}-v^2\eqno(50)$$
in their domains $(-1+\theta,1+\theta)$ and $\big[-\theta,-\theta+(2n-4)^{\un}\big]$. 
\medskip
Observe  that  if $n=4$ then $G(\theta)=B(1-\theta)$.\smallskip

We will be interested in the values $b\ge0$ for which the 
 boundaries $\p B(\theta)$ and $\p G(\theta)$ meet on the hyperplane $y_1=b$. Two special situations occur: when the two boundaries meet on the $y_1$-axis, and when they meet on the hyperplane $y_1=0$; see Figures 3 and 5 in the Appendix. The values of $\theta$ corresponding to those two situations are given as 
$$1+\theta_0=(2n-4)^\un-\theta_0\quad{\hbox{and}}\quad 1=(2n-4)\theta_1^2,$$
that is
$$\theta_0=\theta_0(n)={1\over2}\Big[(2n-4)^\un-1\Big]\qquad{\hbox{and}} \quad\theta_1=\theta_1(n)={1\over \sqrt{2n-4}}.\eqno(51)$$
Note that
$$\cases{\theta_0=\theta_1={1\over2} & for $n=4$\cr
\theta_0>\theta_1 & for $n=5$\cr
\theta_0<\theta_1 & for $n=3$ or $n\ge6.$\cr}$$
\bigskip

In the following  lemmas we will analyze the inclusion relations between the sets $E(\theta),\,D(\theta)$ and the sets $B(\theta),\, G(\theta)$.

\medskip\proclaim Lemma 15.  For any $n\ge 3$ and $\theta\in[0,1]$ we have 
$$E(\theta)\subseteq B(\theta) \qquad D(\theta)\subseteq G(\theta)\eqno(52)$$
\par
\pf Proof of Lemma 15. The result follows instantly from the definition of the four sets.\endpf
 
 \proclaim Lemma 16. For $n\ge 4 $ and  $\theta\ge \theta_1=\ds{1\over\sqrt {2n-4}}$ we have $D(\theta)\subseteq B(\theta)$. If $n=3$ then $D(\theta)=\emptyset$ for $\theta> {1\over2}$.
\par
\pf Proof of Lemma 16. For $n=3$ the condition $|y+\theta e_1|^3\le 2y_1(y_1+\theta)$ implies $(y_1+\theta)^2\le 2y_1$ which has no solutions for $\theta>\half$.

When $n\ge4$ and $y\in D(\theta)$, then $|y+\theta e_1|^2\le (2n-4)^{2/n}\big(y_1(y_1+\theta)\big)^{2/n}$, so that $y\in B(\theta)$ if
$$(2n-4)^{2/n}\big(y_1(y_1+\theta)\big)^{2/n}-4y_1\theta\le1$$
or
$$(1+4y_1\theta)^{n/2}-(2n-4)y_1(y_1+\theta)\ge0$$
under the condition $\theta\ge1/\sqrt{2n-4}$. Since the left hand side is increasing in $\theta$, it is enough to verify
$$\bigg(1+{4y_1\over\sqrt{2n-4}}\bigg)^{n/2}-(2n-4)y_1^2-\sqrt{2n-4}y_1>0$$
or
$$\psi(z):=\Big(1+{2z\over n-2}\Big)^{n/2}-z^2-z>0,\,\qquad z\ge0,\, n\ge4\eqno(53)$$
Note first that 
$$\Big(1+{2z\over n-2}\Big)^{n/2}\ge e^z\qquad z\in[0,2],\;\;n\ge4
$$
indeed it is easy to check that 
$$g(z)={n\over2}\log\Big(1+{2z\over n-2}\Big)-z$$
is concave, $g(0)=0$, and 
$$g(2)={n\over2}\log\Big(1+{4\over n-2}\Big)-2>0.$$
 Therefore, 
$$\psi(z)\ge e^z-z^2-z>0,\qquad 0\le z\le 2,$$
since the function on the right has only one minimum $z_0\in \big[{5\over4},{4\over3}\big],$
and at $z_0$  is greater than  $e^{5\over4}-{16\over9}-{4\over3}>0.$

Now,  for $z\ge2$
$$\psi''(z)={n\over n-2}\Big(1+{2z\over n-2}\Big)^{n/2-2}-2\ge {n\over n-2}\Big(1+{4\over n-2}\Big)^{n/2-2}-2>{n\over n-2}\Big(1+{2(n-4)\over n-2}\Big)-2\ge0$$
so that 
$$\eqalign{\psi'(z)=&{n\over n-2}\Big(1+{2z\over n-2}\Big)^{n/2-1}-2z-1>\psi'(2)={n\over n-2}\Big(1+{4\over n-2}\Big)^{n/2-1}-5\cr&\ge {n\over n-2}\Big(1+{4\over n-2}\Big)^{-1}e^2 -5={n\over n+2} e^2-5>{5e^2\over7}-5>0,\qquad n\ge5\cr}$$
and obviously $\psi'(2)=1>0$ if $n=4$. As a consequence, $\psi(z)\ge \psi(2)>0,$ for any $z\ge 2$.
\endpf
Lemmas 15 and 16 guarantee that for $\theta_1\le \theta\le 1$ both regions  $E(\theta)$ and $D(\theta)$ are  inside~$\,B(\theta)$. The next lemma examines the  relative geometry of $\p B(\theta)$ and $\p G(\theta)$ in more detail; see also Figures 1-5 in the Appendix, which visualize the situation for a generic  $n\ge 6$.\medskip

\proclaim Lemma 17. If $n\ge 4$ and   $0\le \theta<\theta_0$   the boundaries of $B(\theta)$ and $G(\theta)$  intersect on exactly one  hyperplane
 $y_1=b(\theta)$, such that for $\theta>0$
$$0<b(\theta)<\min\Big\{{1\over \theta(2n-4)},1+\theta\Big\}\eqno(54)$$
with the exception $n=5$ and $\theta_1<\theta<\theta_0$, in which case there are no intersections.  Moreover, 
$$B(\theta)\cap\{(y_1,y_2):\,y_1\ge b(\theta)\}\subseteq G(\theta),\eqno(55)$$
with the exception $n=5$ and $\theta_1<\theta<\theta_0$, in which case $B(\theta)\subseteq G(\theta)$. \smallskip
If $n\ge 6$ and   $\;\theta_0\le \theta\le \theta_1$  then the boundaries intersect  on exactly two  hyperplanes  $y_1=b_1(\theta)$ and $y_1=b_2(\theta)$ with 
$$0\le b_1(\theta)<{1\over \theta(2n-4)}<b_2(\theta)\le (2n-4)^\un-\theta\le 1+\theta,\eqno(56)$$ 
 with equality on the left if and only if  $\theta=\theta_1$ and equality on the right if and only if  $\theta=\theta_0$, in which case $b_2(\theta_0)=1+\theta_0$. Moreover, 
$$\eqalignno{&G(\theta)\cap\{(y_1,y_2):\,0\le y_1\le b_1(\theta)\}\subseteq B(\theta)&(57)\cr& 
B(\theta)\cap\{(y_1,y_2):\,b_1(\theta)\le y_1\le b_2(\theta)\}\subseteq G(\theta) &(58)\cr
&G(\theta)\cap\{(y_1,y_2):\,b_2(\theta)\le y_1\le 1+\theta\}\subseteq B(\theta).&(59)\cr}
$$
\par
\pf Proof of Lemma 17. Introduce the function 
$$\phi(b,\theta)=(1+4\theta b)^{n/4}-\sqrt {2n-4}\,(b+\theta),\qquad b,\theta\ge 0.\eqno(60)$$
The intersections between $\partial G(\theta)$ and $\partial B(\theta)$ are given by the equation
$$h(b-\theta)=f(b+\theta)$$
or equivalently $\phi(b,\theta)=0$,
in the range $0\le b\le \min\{1+\theta,(2n-4)^{\un}-\theta\}$. Observe also that 
$$h(b-\theta)>f(b+\theta)\iff \phi(b,\theta)>0.\eqno(61)$$

We already know that $\phi(0,\theta)\ge0$ (i.e. $h(-\theta)\ge f(\theta)$ if $0\le \theta\le \theta_1$, with equality at $\theta_1$, unless $n=5$ and $\theta_1<\theta<\theta_0$, in which case $\phi(0,\theta)<0$.  We also know that if $0\le \theta<\theta_0$ then $1+\theta<(2n-4)^{\un}-\theta$ and $\phi(1+\theta,\theta)<0$ ($h(1)=0< f(1+2\theta)$). Since the function $\phi$ is convex in $b$, this means that there a single zero of $\phi$ on $(0,1+\theta)$ if $\theta\in(0,\theta_0)$, unless $n=5$ and  $\theta_1<\theta<\theta_0$, in which case $\phi(b,\theta)<0$, for $0\le b\le 1+\theta$.

Next, we note that 
$$\phi\Big({1\over \theta(2n-4)},\theta\Big)\le 0,\qquad n\ge 4,\;\; \theta>0,\eqno(62)$$
with equality if and only if $n=4$ and $\theta=\theta_0=\theta_1={\half}$.
Indeed (62) is equivalent  to 
$$ {1\over \theta\sqrt{2n-4}}+\theta\sqrt{2n-4}\ge \bigg(1+{2\over n-2}\bigg)^{n/4}\eqno(63)$$
which is true since  
$$\bigg(1+{2\over n-2}\bigg)^{n/4}\le 2\eqno(64)$$
with equality if and only if $n=4$, in which case equality holds also in (63) precisely when $\theta=\half$. 
 Appliying (62) in the case $n\ge 4$ and $0< \theta<\theta_0$ we obtain that $b(\theta)$ (if it exists) must be also smaller than $1/(\theta(2n-4))$, thereby proving (54).

When  $n\ge 6$ and $\theta_0\le\theta\le \theta_1$ we  have  
$${1\over\theta(2n-4)}< (2n-4)^{\un}-\theta,\qquad n\ge4$$
since it is equivalent to 
 $$ {1\over \theta\sqrt{2n-4}}+\theta\sqrt{2n-4}< (2n-4)^{n\over 2n-4},$$
which in turns follows from  the left hand side being  decreasing  in $\theta\in(0,\theta_1]$ and 
$$\theta\ge \theta_0={\half}\big((2n-4)^{1\over n-2}-1\big)\ge{\log(2n-4)\over 2n-4}\ge{1\over 2n-4}.$$ 
Hence, we conclude that for $n\ge 6$ and $\theta\in[\theta_0,\theta_1]$ we have 
$$\phi(0,\theta)\ge0,\quad\phi\Big({1\over\theta(2n-4)},\theta\Big)<0,\quad \phi\big((2n-4)^\un-\theta,\theta\big)\ge 0$$
(the last inequality being the same as $f\big((2n-4)^\un)=0\le h\big((2n-4)^\un-2\theta\big)$, with equality on the left precisely when $\theta=\theta_1$ and equality on the right when $\theta=\theta_0$. Therefore (56) follows from the convexity of $\phi(\cdot,\theta)$.
\endpf

\proclaim  Lemma 18. The following hold:\medskip
\item {a)} If $n=3$ then
$$E(\theta)\subseteq B\big(\theta e_1,\ts{2\over3}\big),\qquad 0\le\theta\le {1\over2}\eqno(65)$$
\medskip
\item{b)} If $n\ge4$ then 
$$ E(\theta)\subseteq G(\theta),\qquad 0\le\theta<\theta_0\eqno(66)$$
and
$$E(\theta)\cap\{y:\,y_1\theta(2n-4)\le 1\}\subseteq G(\theta),\qquad 0\le\theta\le 1.\eqno(67)$$

 \par

\pf Proof of Lemma 18.  We begin by noting  that if $y\in E(\theta)$ then
$$|y-\theta e_1|^{2-n}-|y+\theta e_1|^{2-n}\ge 1\eqno(68)$$
and $|y-\theta e_1|\le 1$.
This last estimate can be improved a tad as follows:
$$|y+\theta e_1|^2=|y-\theta e_1|^2+4y_1\theta\le 1+4(1+\theta)\theta=(1+2\theta)^2$$
so that $$|y-\theta e_1|^{2-n}\ge 1+(1+2\theta)^{2-n}$$
i.e. $$|y-\theta e_1|\le \Big(1+(1+2\theta)^{2-n}\Big)^{-\un}.$$
Using this estimate for $n=3$ and $\theta\le\half$ gives (65).

To show (66), we  start with a preliminary inclusion. Define 
$$B^*(\theta)=\{y\in\Rn:\, |y-\theta e_1|^n\le (2n-4)\theta y_1\}$$
and let us show that 
$$E(\theta)\subseteq B^*(\theta),\qquad \theta\ge0.\eqno(69)$$
From (68) we get 
$$|y+\theta e_1|^{n-2}(1-|y-\theta e_1|^{n-2})\ge|y-\theta e_1|^{n-2}$$
if $R=|y-\theta e_1|<1$ then the above inequality  becomes
$$(R^2+4 y_1\theta)^{n-2\over2}\ge{R^{n-2}\over 1-R^{n-2}}$$
$$4 y_1\theta\ge\bigg({R^{n-2}\over 1-R^{n-2}}\bigg)^{2\over n-2}-R^2=R^2\bigg[\bigg({1\over 1-R^{n-2}}\bigg)^{2\over n-2}-1\bigg]\ge {2\over n-2} R^n,$$
which is (69).

At this point we know that 
$$E(\theta)\subseteq B^*(\theta)\cap B(\theta,1)=\Big\{y:|y-\theta e_1|^n\le\min\{1,(2n-4)y_1\theta\}\,\Big\}.$$
 
Now for a point $y\in E(\theta)$, for any $\theta\in [0,1]$, we have either
$$ y_1\theta (2n-4)\le 1\quad{\hbox{ and }}\quad  y\in B^*(\theta)\eqno(70)$$
or
$${1\over \theta(2n-4)}< y_1\le 1+\theta\quad{\hbox{ and }}\quad  y\in B(\theta),\eqno(71)$$
in the assumption that  $1+\theta> 1/\big( \theta(2n-4)\big)$, i.e.
$\theta> \theta_{00}:={1\over2}\Big[\Big(1+{2\over n-2}\Big)^{1/2}-1\Big]$ (which is when the surfaces $\partial B^*(\theta)$ and $\partial B(\theta)$ intersect at    $y_1=1/(\theta(2n-4))$.)

If (70) holds for some $\theta\in[0,1]$, then  
$$|y+\theta e_1|^2=|y-\theta e_1|^2+4\theta y_1\le (2n-4)^{2/n}(y_1\theta)^{2/n}+4y_1\theta,$$
and  $y\in G(\theta)$ provided
$$(2n-4)^{2/n}(y_1\theta)^{2/n}+4y_1\theta\le (2n-4)^{2/n}(y_1+\theta)^{4/n}.$$
But $$\eqalign{(2n-4)^{2/n}\Big((y_1+\theta)^{4/n}-(y_1\theta)^{2/n}\Big)&= (2n-4)^{2/n}\Big((y_1^2+2y_1\theta+\theta^2)^{2/n}-(y_1\theta)^{2/n}\Big)\cr &\ge(2n-4)^{2/n}(4^{2/n}-1) (y_1\theta)^{2/n}\ge 4 y_1\theta\cr}$$
provided 
$$y_1\theta\le (2n-4)^{2\over n-2}\bigg({4^{2/n}-1\over 4}\bigg)^{n\over n-2}.$$
However   
$$y_1\theta\le{1\over 2n-4}\le(2n-4)^{2\over n-2}\bigg({4^{2/n}-1\over 4}\bigg)^{n\over n-2} $$ 
as the last inequality is equivalent to (64). This shows that a point $y\in E(\theta)$ is also in $G(\theta)$  in case (70) holds, thereby proving (67).

If instead $y\in E(\theta)$, $\theta_{00}<\theta< \theta_0$ and (71) holds, then (54) and (55) immediately imply that $y\in G(\theta)$, and this, together with (67), proves (66).

\smallskip

 We can now summarize the results obtained in Lemmas 15, 16 and 18  in the following: 
\bigskip\proclaim Corollary 19. If $\theta_0$ and $\theta_1$ are as in (51), then \bigskip\item{a)} If $n=3$
$$E(\theta)\cup D(\theta)\subseteq\cases{ B\big(\theta e_1,{2\over3}\big)\cup G(\theta) & if  $0\le \theta\le{1\over2}$\cr
 B(\theta) & if  ${1\over2}\le \theta\le 1.$\cr}$$
$$|E(\theta)\cup D(\theta)|\le |B_3|,\qquad 0\le\theta\le 1.\eqno(72)$$
\smallskip\item{b)} If $n=4$ or $n=5$
 $$E(\theta)\cup D(\theta)\subseteq\cases{  G(\theta) & if  $0\le \theta<\theta_0$\cr
 B(\theta) & if  $\theta_1\le \theta\le 1$.\cr}$$
$$|E(\theta)\cup D(\theta)|\le \cases{|B_n| & if $n=4$\cr|G_n| & if $n=5$\cr},\qquad 0\le\theta\le 1$$
\smallskip
\item{c)}  If $n\ge6$
 $$E(\theta)\cup D(\theta)\subseteq\cases{  G(\theta) & if  $0\le \theta<\theta_0<\theta_1$\cr
 B(\theta) & if  $\theta_1\le \theta\le 1$.\cr}$$
$$|E(\theta)\cup D(\theta)|\le |G_n|,\qquad \theta\in[0,\theta_0]\cup[\theta_1,1].$$
\par
\pf Proof of Corollary 19. The only thing to check here is (72):
$$|E(\theta)\cup D(\theta)|\le |B_3| \bigg({2\over3}\bigg)^3+|G_3|={4\pi\over3}\,{8\over27}+{16\pi\over 21}={4\pi\over3}\,{164\over 189}<{4\pi\over3}=|B_3|.$$\endpf

\bigskip
It is clear from the previous corollary that the  only gap remaining 
toward a complete proof of Proposition 14, is the volume estimate in the case $\theta_0\le\theta\le \theta_1$ and $n\ge6$.
Numerical evidence shows that in that range of $\theta$'s
and for large enough $n$, it is in general false that  $E(\theta)$ 
and $D_0(\theta)$ (a proper subset of $D(\theta)$) are either both inside $B(\theta)$ or both inside $G(\theta)$. Thus it seems hopeless to try to play with inclusions in order to give an estimate for $|E(\theta)\cup D_0(\theta)|$. Nonetheless, we are able to show what we need:

\proclaim Lemma 20. If $n\ge6$ and $\theta_0\le\theta\le\theta_1$ we have
$$|E(\theta)\cup D(\theta)|\le |G_n|.$$
\par\medskip

\pf Proof of Lemma 20. We know from Lemma 17 that the equation $\phi(b,\theta)=0$  that gives the intersection of $\p B(\theta)$ and $\p G(\theta)$ has two distinct solutions $b_1=b_1(\theta)$ and $b_2=b_2(\theta)$ as in (56) (see Fig. 3,4,5). From (56) and (67) we have 
$$\{y:E(\theta),\, y_1\in[0,b_1]\}\subseteq G(\theta),$$
from (58) we have
$$\{y\in E(\theta):\, y_1\in[b_1,b_2]\}\subseteq \{y\in B(\theta):\, y_1\in[b_1,b_2]\}\subseteq G(\theta)$$
and clearly
$$\{y:\,y\in E(\theta),\, y_1\in[b_2,1+\theta]\}\subseteq B(\theta)$$
 $$\{y:\,\,y\in D(\theta),\,y_1\in\big[0,b_2]\}\subseteq  G(\theta).$$
Finally, from (59) 
$$\{y:\,\,y\in D(\theta),\,y_1\in\big[b_2,1+\theta]\}\subseteq  B(\theta).$$

This means that we can use the following volume bound:
$$\eqalign{|E(\theta)\cup D(\theta)|&\le |\{y\in G(\theta):\, y_1\in[0,b_2]\}|+|\{y\in B(\theta):\, y_1\in[b_2,1+\theta]\}|\cr& ={\omega_{n-2}\over n-1}\int_0^{b_2} f(y_1+\theta)^{n-1\over2}dy_1+{\omega_{n-2}\over n-1}\int_{b_2}^{1+\theta} h(y_1-\theta)^{n-1\over2}dy_1\cr}$$
that is 
$$|E(\theta)\cup D(\theta)|\le V(\theta)$$
where 
$$V(\theta)={\omega_{n-2}\over n-1}\int_\theta^{b_2+\theta} f(v)^{n-1\over2}dv+{\omega_{n-2}\over n-1}\int_{b_2-\theta}^{1} h(v)^{n-1\over2}dv$$
where $f$ and $h$ are defined in (50).

The goal is to show that $V(\theta)\le |G(0)|$. This inequality is obvious at $\theta=\theta_0$, since at $\theta_0$ the second integral vanishes ($b_2(\theta_0)=1+\theta_0$), so it would be enough to show that $V$ is decreasing on $[\theta_0,\theta_1]$, but unfortunately this fact turns out to  be true only for $n\le 12$. What we show instead is that $V$ has at most one extremum, which  is a minimum, and that $V(\theta_1)\le |G(0)|$.
\eject
We have
$${n-1\over \omega_{n-2}}\,V'(\theta)=(b_2'+1) f(b_2+\theta)^{n-1\over2}-f(\theta)^{n-1\over2}-(b_2'-1)h(b_2-\theta)^{n-1\over2}=2h(b_2-\theta)^{n-1\over2}-f(\theta)^{n-1\over2}$$

and (note that $b_2(\theta)>\theta$, due to (56))
$$V'(\theta)<0\Longleftrightarrow 2^{2\over n-1}\Big(1-(b_2-\theta)^2\Big)<f(\theta)\Longleftrightarrow b_2(\theta)>q(\theta),$$
where 
$$q(\theta)=\theta+\sqrt{1-2^{-{2\over n-1}} f(\theta)}.$$
Taking into account Lemma 17 and (61)
$$b_1(\theta)<q(\theta)<b_2(\theta)\Longleftrightarrow \phi(\theta):=\phi(q(\theta),\theta)<0$$
where $\phi(b,\theta)$ is defined in (60).

We now show that $\phi(\theta)$ is strictly increasing, so it has at most one zero. We  prove
$$\phi'(\theta)=n(q+\theta q')(1+4\theta q)^{n/4-1}-\sqrt{2n-4}\, (q'+1)>0,\qquad \theta_0\le\theta\le\theta_1\eqno(73)$$
where $q$ and $q'$ are evaluated at $\theta$. 

\proclaim Claim. For $\theta_0\le \theta\le\theta_1$ we have 
$$\hskip3.3em q'(\theta)<0,\qquad n\ge 9\eqno(74)$$
$$q(\theta)+\theta q'(\theta)>0,\qquad n\ge 6.\eqno(75)$$ 
\par\medskip

Assuming the above claim,  for  $n\ge 9$  and    $\theta_0\le \theta\le \theta_1$ we have 
$$\eqalign{\phi'(\theta)&>n(q+\theta q')-\sqrt{2n-4}\,(q'+1)=
n q-\sqrt{2n-4}+ (n\theta-\sqrt{2n-4})q'\cr&>n q(\theta_1)-\sqrt{2n-4}=n\sqrt{1-2^{-{2\over n-1}}{2n-5\over 2n-4}}-{n-4\over \sqrt{2 n-4}}.\cr}$$

So $\phi'(\theta)>0$ if 
$$n^2\big(2n-4-2^{-{2\over n-1}}\,(2n-5)\big)>(n-4)^2$$
or 
$$2n^3\big(1-2^{-{2\over n-1}}\big)-5n^2(\big(1-2^{-{2\over n-1}}\big)+8n-16>0$$
or
$$\big(1-2^{-{2\over n-1}}\big)n^2(2n-5)+8n-16>0,$$
which is obvious for $n\ge4$. This settles (73) when $n\ge 9$.

To deal with the cases $n=6,7,8$ (and in those cases it's not true that $q'<0$ on $[\theta_0,\theta_1]$),
start by writing (again with $\theta_0\le \theta\le \theta_1$)
$$\phi'(\theta)>n(q+\theta q')-\sqrt{2n-4}\,(q'+1)>n q+n\theta-2\sqrt{2n-4}$$
since $q'<1$, and $n\theta-\sqrt{2n-4}<0$. All we need to show is that $nq>2\sqrt{2n-4}-n\theta$ or 
$$f(\theta)<1-2^{2\over n-1}\Big({2\over n}\sqrt{2n-4}-2\theta\Big)^2,\qquad \theta_0\le\theta\le  \theta_1$$
for $n=6,7,8$. This is implied by
$$f(\theta_1)\le 1-2^{2\over n-1}\Big({2\over n}\sqrt{2n-4}-2\theta_0\Big)^2,\qquad n=6,7,8$$
which can be verified numerically.

This shows (assuming the Claim)  that $\phi$ is strictly increasing in $[\theta_0,\theta_1]$. Now observe that the function $f$ has a maximum at $\theta=\big({2\over n}\big)^{n\over 2n-4}(2n-4)^{1\over n-2}>\theta_1$ and $f(\theta_1)<1$.
This means that in the range $\theta_0\le\theta\le \theta_1$ we have $\theta<q(\theta)<1+\theta$, and in particular 
$b_2(\theta_0)=1+\theta_0>q(\theta_0)$. We claim that $\phi(\theta_0)<0$. If it were $\phi(\theta_0)>0$, then $\phi(\theta)>0$ for all $\theta\in [\theta_0,\theta_1]$, which implies that $0\le q(\theta)<b_1(\theta)$ for all such $\theta$'s. But that is not possible since it would imply  $q(\theta_1)=b_1(\theta_1)=0$, by continuity of  $b_1$. 

Since $\phi(\theta_0)<0$ then $\phi$ is negative on $[\theta_0,\theta_1]$ provided  $\phi(\theta_1)<0$, and this can be checked numerically  if $6\le n\le 12$. For $n\ge 13$ one could prove  that $\phi(\theta_1)>0$, however this is not necessary for our purposes (the reader can verify for example that $\phi(\theta_1)\to e^{\half\sqrt{1+\log16}}-\sqrt{1+\log16}-2>0$, as $n\to +\infty$).   Indeed, we know that since $\phi(\theta)$ has at most one zero, and it's negative at $\theta_0$,  then   $V$ has at most  one minimum on $[\theta_0,\theta_1]$ if $n\ge 13$ and it's decreasing in that interval for $n\le12$.  Since $V(\theta_0)\le|G(0)|$ it is now enough to prove that  $V(\theta_1)\le |G(0)|$ for $n\ge 13.$ The inequality is written as

$$\int_{\theta_1}^{b_2+\theta_1}
f(v)^{n-1\over2}dv+\int_{b_2-\theta_1}^1
 h(v)^{n-1\over2}dv\le \int_0^{(2n-4)^{1\over n-2}} f(v)^{n-1\over2}dv\eqno(76)$$
where $b_2=b_2(\theta_1)$ is the only positive solution of the equation
$$(1+4\theta_1 b)^{n/4}-{b\over\theta_1}-1=0.$$

Make the change $v=\theta_1 x$ and obtain that (76) is equivalent to
$$\eqalign{ \int_1^{\lambda_n+1} \bigg(x^{4/n}-{x^2\over 2n-4}\bigg)^{n-1\over2}dx+&\int_{\lambda_n-1}^{\sqrt {2n-4}}\bigg(1-{x^2\over 2n-4}\bigg)^{n-1\over2}dx\le\cr&\le \int_0^{(2n-4)^{{1\over2}+{1\over n-2}}}\bigg(x^{4/n}-{x^2\over 2n-4}\bigg)^{n-1\over2}dx\cr}\eqno(77)$$ 
where $\lambda_n$ is the unique positive solution of the equation
$$\bigg(1+{2\lambda\over n-2}\bigg)^{n/4}-1-\lambda=0.
$$
Rewrite (77) as
$$\eqalign{J(n):=\int_0^1 x^{2-2/n}\bigg(1-{x^{2-4/n}\over 2n-4}\bigg)^{n-1\over2}dx&+\int_{\lambda_n+1}^{(2n-4)^{{1\over2}+{1\over n-2}}} x^{2-2/n}\bigg(1-{x^{2-4/n}\over 2n-4}\bigg)^{n-1\over2}dx-\cr&\hskip 4em-\int_{\lambda_n-1}^{\sqrt {2n-4}}\bigg(1-{x^2\over 2n-4}\bigg)^{n-1\over2}dx\ge0\cr}$$
First notice that  if $w=A/(2n-4)$ then
$${\partial\over\partial n}(1-w)^{n-1\over2}={1\over2}(1-w)^{n-1\over2}\bigg({n-1\over n-2}\cdot{w\over 1-w}+\log(1-w)\bigg)>0\eqno(78)$$
for $0<A<2n-4$, and $n>2$. Next, if 
$$g(\lambda,n)= \bigg(1+{2\lambda\over n-2}\bigg)^{n/4}-1-\lambda$$
then it is straightforward to check that $g$ is increasing in $n$ for $n\ge 13$ and $\lambda>2.5.$
In particular, 
$$g(\lambda,n)\le e^{\lambda/2}-1-\lambda:=g(\lambda)\qquad n\ge 13.$$

\medskip
These last facts allow us to localize the values $\lambda_n$:

$$g(2.51, n)\le g(2.51)\approx -0.0021,\qquad n\ge 13$$
$$g(2.56, n)\ge g(2.56,66)\approx 0.00034,\qquad n\ge 66$$
$$g(2.56, n)\le g(2.56, 65)\approx -0.00017,\qquad n\le 65$$
$$g(2.61, n)\ge g(2.61, 33)\approx 0.00069,\qquad n\ge 33$$
$$g(2.61, n)\le g(2.61, 32)\approx -0.0013,\qquad n\le 32$$
$$g(2.67, n)\ge g(2.67, 21)\approx 0.00051,\qquad n\ge 21$$
$$g(2.67, n)\le g(2.67, 20)\approx -0.0044,\qquad n\le 20$$
$$g(2.79,n)\ge g(2.79,13)\approx 0.0042,\qquad n\ge 13$$

and these relations imply\eject
$$\hskip-2.3em 2.51<\lambda_n<2.56,\qquad n\ge 66$$
$$2.56<\lambda_n<2.61,\qquad 33\le n\le 65$$
$$2.61<\lambda_n<2.67,\qquad 21\le n\le 32$$
$$2.67<\lambda_n<2.79,\qquad 13\le n\le 20.$$
Now, if $n\ge n_1=66$ we have $x^{2-2/n}\ge x^2$ and $x^{2-4/n}\le x^{2-4/n_1}$, for $0<x<1$, and $x^{2-2/n}\ge x^{2-2/n_1},\;$ $x^{2-4/n}\le x^2$ for $x>1$, hence, taking into account (78), 
$$\eqalign{ J(n)\ge\int_0^1 x^2 \bigg(1-{x^{2-4/n_1}\over 2n_1-4}\bigg)^{n_1-1\over2}dx+\int_{3.56}^{\sqrt{2n_1-4}} x^{2-2/n_1}&\bigg(1-{x^{2}\over 2n_1-4}\bigg)^{n_1-1\over2}dx-\cr&-\int_{1.51}^\infty e^{-x^2/4}dx\approx 0.0018.\cr}$$
For any $n_2>n_1\ge 13$, if  $n_1\le n\le n_2$ and $\mu_1<\lambda_n<\mu_2$ we have
$$\eqalign{ J(n)\ge\int_0^1 x^{2-2/n_2} \bigg(1-{x^{2-4/n_1}\over 2n_1-4}\bigg)^{n_1-1\over2}dx&+\int_{\mu_2+1}^{\sqrt{2n_1-4}} x^{2-2/n_1}\bigg(1-{x^{2-4/n_2}\over 2n_1-4}\bigg)^{n_1-1\over2}dx-\cr&\qquad -\int_{\mu_1-1}^{\sqrt{2n_2-4}} \bigg(1-{x^2\over 2n_2-4}\bigg)^{n_2-1\over2}dx.\cr}$$
Using this estimate and the above bounds on $\lambda_n$  we find
$$J(n)\ge \cases{0.030 & if $33\le n\le 65$\cr
0.046 & if $21\le n\le 32$\cr
0.018 & if $13\le n\le 20$\cr}$$
and this shows that $J(n)>0$ for $n\ge 13$, concluding the proof of Lemma 20.
\endpf
 Lemma 20 concludes the proof of Proposition 14, and hence the proofs of estimates (40), (41) and Theorem 7.\endpf
\bigskip
\eject
\pf Proof of Claim. 
\medskip
 We begin by proving (74). Let $F(\theta)=2^{-{2\over n-1}} f(\theta)$, so that $q(\theta)=\theta+\sqrt {1-F(\theta)}$ and
$$\qquad q'<0\Longleftrightarrow F'>2\sqrt{1-F}\Longleftrightarrow (F')^2>4(1-F),$$
since $F$ is increasing in our range. This last estimate is proven once we show that 
$$ {d\over d\theta} \Big[(F')^2-4(1-F)\Big]=2F'(F''+2)<0,\quad n\ge5\eqno(79)$$
and
$$F'(\theta_1)>2\sqrt{1-F(\theta_1)},\quad  n\ge9.\eqno(80)$$

 To show (79), i.e. $F''+2<0$ on $[\theta_0,\theta_1]$, note that 
$$F'(\theta)=2^{-{2\over n-1}}\Big[(2n-4)^{2\over n}{4\over n} \theta^{4/n-1}-2\theta\Big]>0,\qquad \theta\in[\theta_0,\theta_1]$$
$$F''(\theta)=2^{-{2\over n-1}}\Big[(2n-4)^{2\over n}{4\over n}\bigg({4\over n}-1\bigg) \theta^{4/n-2}-2\Big]<0,\qquad \theta\in[\theta_0,\theta_1]$$
$$F'''(\theta)=2^{-{2\over n-1}}(2n-4)^{2\over n}{4\over n}\bigg({4\over n}-1\bigg)\bigg({4\over n}-2\bigg) \theta^{4/n-3}>0,\qquad \theta\in[\theta_0,\theta_1]$$
so that  
$$\eqalign{F''(\theta)+2&<F''(\theta_1)+2=-2^{-{2\over n-1}}\Big[(2n-4)^{2\over n}{4\over n}\bigg({n-4\over n}\bigg) (2n-4)^{-{2\over n}+1}+2\Big]+2\cr&=-2^{-{2\over n-1}+1}\Big[{4(n-2)(n-4)\over n^2}+1\Big]+2\cr}$$
 We then only need to check whether
$$-2^{-{2\over n-1}}\Big[{4(n-2)(n-4)\over n^2}+1\Big]+1<0,\qquad n\ge5$$
or
$$\Big(5-2^{2\over n-1}\Big)n^2-24n+32>0,\quad n\ge5,$$
which is easy to do.

Estimate (80), after squaring and simplifying,  is equivalent to
$$2n^3(2^{2\over n-1}-2^{4\over n-1})+n^2(9+ 4\cdot  2^{4\over n-1}-5\cdot  2^{2\over n-1})-48 n +64>0,\qquad n\ge9.$$
This inequality can be checked directly for $n=9,10$, and for $n\ge 11$ we can easily argue as follows.
The coefficient of $n^2$ is greater than 8, while for some $u^*\in\big({2\over n-1},{4\over n-1})$
$$2^{2\over n-1}-2^{4\over n-1}=-{2\over n-1}\; 2^{u^*} \log2 >-{2\over n-1}\, 2^{4\over n-1} \log2>-{2^{7/5} \log 2\over n-1},\qquad n\ge11, $$
so we are reduced to check whether
$$-2{2^{7/5} \log 2\over n-1}\,n^3+8n^2-48n+64>0,\qquad n\ge11$$
which is implied by 
$$\Big(-{11\over5}\,2^{7/5}\log 2 +8\Big)n^2-48 n+64>0,\qquad n\ge 11,$$
and this shows (74).

To prove (75) write  
$$q+\theta q'=\theta+\sqrt{1-F}+\theta\bigg(1-{F'\over2\sqrt{1-F}}\bigg)>0\;\iff\; 4\theta\sqrt{1-F}+2(1-F)-\theta F'>0$$
so it's enough to prove 
$$2(1-F)-\theta F'>0\qquad \theta\in[\theta_0,\theta_1],\eqno(81)$$
and in particular it's enough to prove that the left-hand side is decreasing on $[\theta_0,\theta_1]$, and that the inequality above is verified at $\theta=\theta_1$.

It's easy to check that the derivative of the left-hand 
side of (81) coincides with
$$-3F'-\theta F''=-8\theta 2^{-{2\over n-1}}\bigg({n+2\over n^2}\,(2n-4)^{2/n}\theta^{4/n-2}-1\bigg)$$
so $-3F'-\theta F''<0$ if 
$$\theta<\bigg({n+2\over n^2}\bigg)^{n\over 2n-4}(2n-4)^{1\over n-2}\qquad\theta_0\le\theta\le \theta_1$$
but this condition is easily verified if $n\ge 6$ since the right-hand side of the above inequality is larger than $\theta_1$. 

Now we only need to check that (81) holds at $\theta=\theta_1$, but this is easy since
$$2\big(1-F(\theta_1)\big)-\theta_1 F'(\theta_1)=2\bigg[1+2^{-{2\over n-1}}\bigg({1\over n-2}-{2\over n}-1\bigg)\bigg]\ge 2^{1-{2\over n-1}} \;{2\log 2-1\over n}>0,$$
and this conludes the proof of $q+\theta q'>0$ for $n\ge 6$.
\endpf
\centerline{\bf 5. Conclusion: Proofs of Theorems 1,2,6}
\medskip
 
As we noted earlier Theorem 1 is a special case of Theorem 2, and Theorem 2 follows from Theorem 6 and Theorem 5,  since the operator $T$ is a bijection. The inequality statement of Theorem 6, on the other hand, is a consequence of Theorem 4 and the distribution estimates of  Theorem 7. The only thing left to prove is the sharpness statement of Theorem 6. In order to do that, we apply the sharpness result of Theorem 4:  equality in (21) and (22) is attained at any  $x_0\in B_n$ for $n=3,4$ and at $x_0\in \p B_n$, for $n\ge 4$. 
We only treat the case  $n\ge4$, under the hypothesis that  there exists $x_0\in \p B_n$ such that $\nu(B(x_0,r)\cap\bar B_n)\ge C_1 r^\lambda$, for $0<r\le r_1$, some $C_1,r_1>0$. The argument  for $n=3,4$,  with the above condition on $\nu$ verified for $x_0\in B_n$,   is similar, and easier.

We can assume   $x_0=e_1=(1,0,...,0)$,  we take   $m$ large enough so that 
$$ \{z\in B_n:\,|K(e_1,z)|>m\,\}\subseteq B(e_1, Cm^{-p'/n})\cap \B_n,$$
and we  let $$r_m=Cm^{-p'/n},\;\; E_m=B(e_1,r_m)\cap\B_n\;\;F_m=B\big(e_1,\ts{1\over10}r_m\big)\cap \B_n.$$
 Conditions (15),(16), (17) of Theorem 4 are met, with $\beta=n/(n- d)$ and $\beta_0=\lambda/(n- d)$, given the hypothesis on $\nu$, so all we need to check is (18), i.e. we will prove the following H\"ormander type condition

$$\int_{|z|\le1, |z-e_1|\ge r_m} |K(x,z)-K(e_1,z)|\,|K(e_1,z)|^{2/(n-2)}dz\le C\eqno(82)$$
for all $x\in\B_n$ with $|x-e_1|<r_m/10$. Given the asymptotic estimate (31) it will suffice to prove (82) for 
$$K(x,z)=c_n|x-z|^{2-n}-c_n\big(1+g(x,z)\big)\big|x^*-|x|z\big|^{2-n}$$ where
$$g(x,z)=2(n-2)\,{x^*\cdot(x^*-z)\,x^*\cdot (x^*-|x|z)\over |x^*-|x|z|^2}.$$
Estimate (82) is a consequence of the following:
$$\int_{|z|\le1, |z-e_1|\ge r_m}\Big||x-z|^{2-n}-|e_1-z|^{2-n}\Big|\,|e_1-z|^{-2}dx\le C\eqno(83)$$
 $$\int_{|z|\le1, |z-e_1|\ge r_m}\Big||x^*-|x|z|^{2-n}-|e_1-z|^{2-n}\Big|\,|e_1-z|^{-2}dx\le C\eqno(84)$$
$$\int_{|z|\le1, |z-e_1|\ge r_m}|g(x,z)-g(e_1,z)|\,|e_1-z|^{-n}dz\le C\eqno(85)$$
for  $|e_1-x|\le r_m/10$.

Inequality (83) is derived  using the estimate
$$\Big||x-z|^{2-n}-|e_1-z|^{2-n}\Big|\,|e_1-z|^{-2}\le C|x-e_1| |e_1-z|^{-n-1}\eqno(86)$$
which is valid under our assumptions, and more generally if $|x-e_1|\le \delta r_m$ any $\delta<1$.

Inequality (84) is a consequence of the estimate
$$\Big||x^*-|x|z|^{2-n}-|e_1-z|^{2-n}\Big|\,|e_1-z|^{-2}\le C|x-e_1| |e_1-z|^{-n-1}$$
which can be asily derived from (86):
$$\eqalign{\Big||x^*-|x|z|^{2-n}-|e_1-z|^{2-n}\Big|\le |x|^{n-2}\bigg|\Big|{x^*\over |x|}-z\Big|^{2-n}-|e_1-z|^{2-n}\Big|+|e_1-z|^{2-n}\Big|1-|x|^{n-2}\Big|\cr}$$
now for large $m$ we have $1/2<|x|\le 1$, and
$$\Big|e_1-{x^*\over |x|}\Big|\le \Big|e_1-x+x-{x\over|x|^2}\Big|\le |e_1-x|+{1-|x|^2\over|x|}\le |e_1-x|\Big(1+{|e_1+x|\over |x|}\Big)\le 5|e_1-x|<{r_m\over2}$$
so that inequality (86) applies with $x^*/|x|$ in place of $x$. Note also that $\Big|1-|x|^{n-2}\Big|\le C|e_1-x|$.

We now only  need to check (85). The numerator of $g(x,z)-g(e_1,z)$ is equal to 
$$\eqalign{&|e_1-z|^2  x^*\cdot(x^*-z) x^*\cdot(x^*-|x|z)-|x^*-|x|z|^2 \big(e_1\cdot(e_1-z)\big)^2=\cr&=
\big(e_1\cdot(e_1-z)\big)^2\Big[|e_1-z|^2-|x^*-|x|z|^2\Big] +\cr & \hskip5em +|e_1-z|^2\Big[x^*\cdot(x^*-z) x^*\cdot(x^*-|x|z)-\big(e_1\cdot(e_1-z)\big)^2\Big]\cr
}\eqno(87)
$$

Now,
$$\eqalign{|e_1-z|^2-|x^*-|x|z|^2& =2x\cdot z-2z\cdot e+|z|^2(1-|x|^2)\cr
&=2z\cdot(x-e_1)+|z|^2(e_1-x)\cdot(e_1+x)=(x-e_1)\cdot\big[2z-|z|^2(x+e_1)\big]\cr&=
(x-e_1)\cdot\big[2(z-e_1)-(x- e_1)|z|^2-2e_1(z-e_1)\cdot(z+e_1)\big]
}
$$
so the first term of (87) is bounded above by 
$$C|z-e_1|^2|x-e_1|^2+C|z-e_1|^3|x-e_1|\eqno(88)$$
For the second term we have
$$\eqalign{x^*\cdot(x^*-z)& x^*\cdot(x^*-|x|z)-\big(e_1\cdot(e_1-z)\big)^2=\cr &=1-x^*\cdot z-x\cdot z+(x^*\cdot z)(x\cdot z)-1+2e_1\cdot z+(e_1\cdot z)^2\cr&=-z\cdot(x-e_1)-z\cdot(x^*-e_1)+(x^*\cdot z)z\cdot (x-e_1)+(e_1\cdot z)z\cdot(x^*-e_1)\cr&= z\cdot(x-e_1)\big[(x^*\cdot z)-1\big]+z\cdot(x^*-e_1)\big[(e_1\cdot z)-1\big]\cr&= z\cdot(x-e_1)\big[(x^*-e_1)\cdot z+ (z-e_1)\cdot e_1\big]+z\cdot(x^*-e_1)\big[(e_1\cdot (z-e_1)\big]}$$
Noting that for $m$ large
$$|x^*-e_1|\le 2\big|x-|x|e_1|\big|\le 2|x-e_1|+1-|x|^2\le 4|x-e_1|$$ we get that 
the second term in (87) is also  bounded above by the quantity in (88). In summary, 
$$|g(x,z)-g(e_1,z)|\le C{|x-e_1|^2|z-e_1| ^2+|z-e_1|^3|x-e_1|\over |x^*-|x|z|^2|z-e_1|^2}$$
and since $|x^*-|x|z|\ge |x-z|\ge|z-e_1|-|x-e_1|$ we have
$$|g(x,z)-g(e_1,z)|\,|z-e_1|^{-n}\le C{|x-e_1|^{2}|z-e_1|^{-n}+|x-e_1||z-e_1|^{-n+1}\over(|z-e_1|-|x-e_1|)^2} $$
and it's now easy to check that (85) holds.
This concludes the proof of the sharpness statement, and hence the proof of Theorem~6.\endpf
\smallskip
\centerline{\bf Appendix}
\medskip
We present a few graphs of the boundaries of $G(\theta)=\overline G_n+\theta e_1$ and $B(\theta)=\overline B_n-\theta e_1$, restricted to the 2-dimensional quadrant $\{y_1\ge0,\,y_2\ge0\}$, for some critical ranges  of~$\theta$. These graphs were plotted with Mathematica when $n=15$, but the pattern is similar for any $n\ge 6$. The dotted line represents $\p B(\theta)$ and the continuous line represents $\p G(\theta)$. The notation for the coordinates where the boundaries intersect is the same as that of Lemma~17.

\input insbox
\input epsf

\def\box1{\vbox{
  \epsfxsize=9cm
  \hbox{\epsfbox{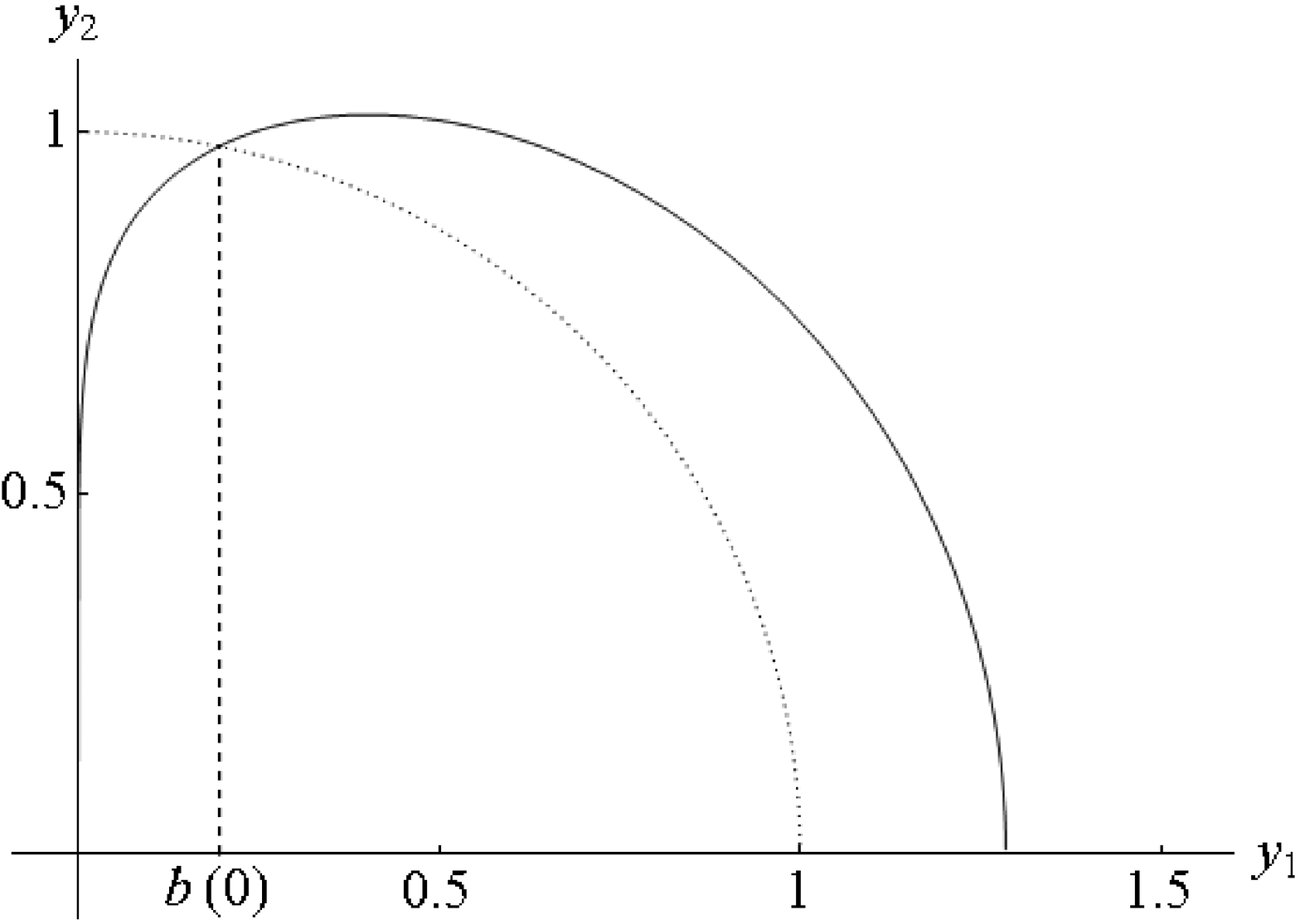}}
  \hbox to 9cm{\hss Figure 1:  $\theta=0$.\hss}
}}

\InsertBoxC{\box1}

\def\box2{\vbox{
  \epsfxsize=9cm
  \hbox{\epsfbox{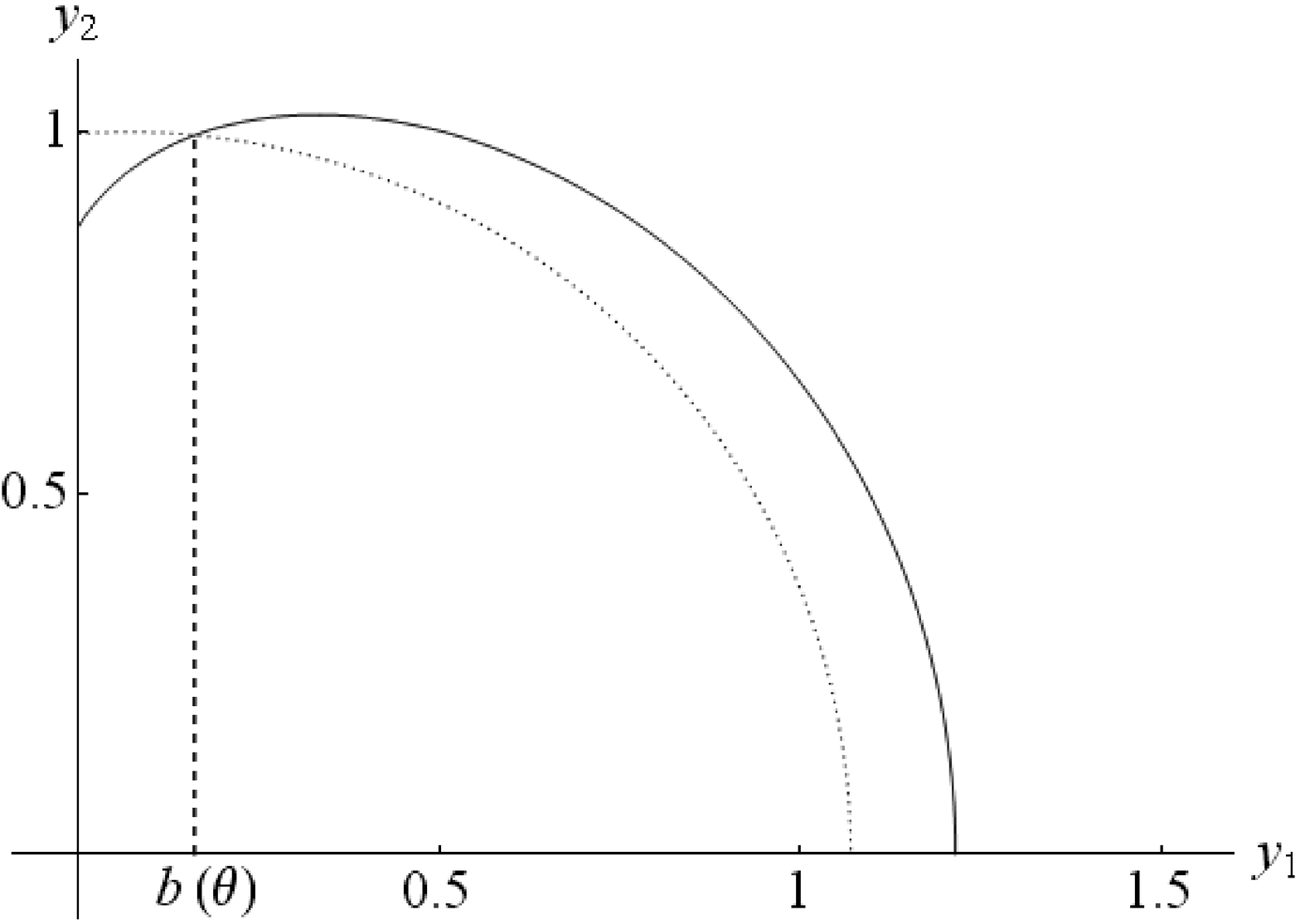}}
  \hbox to 9cm{\hss Figure 2:  $0<\theta<\theta_0$.\hss}
}}

\InsertBoxC{\box2}

\def\box3{\vbox{
  \epsfxsize=9cm
  \hbox{\epsfbox{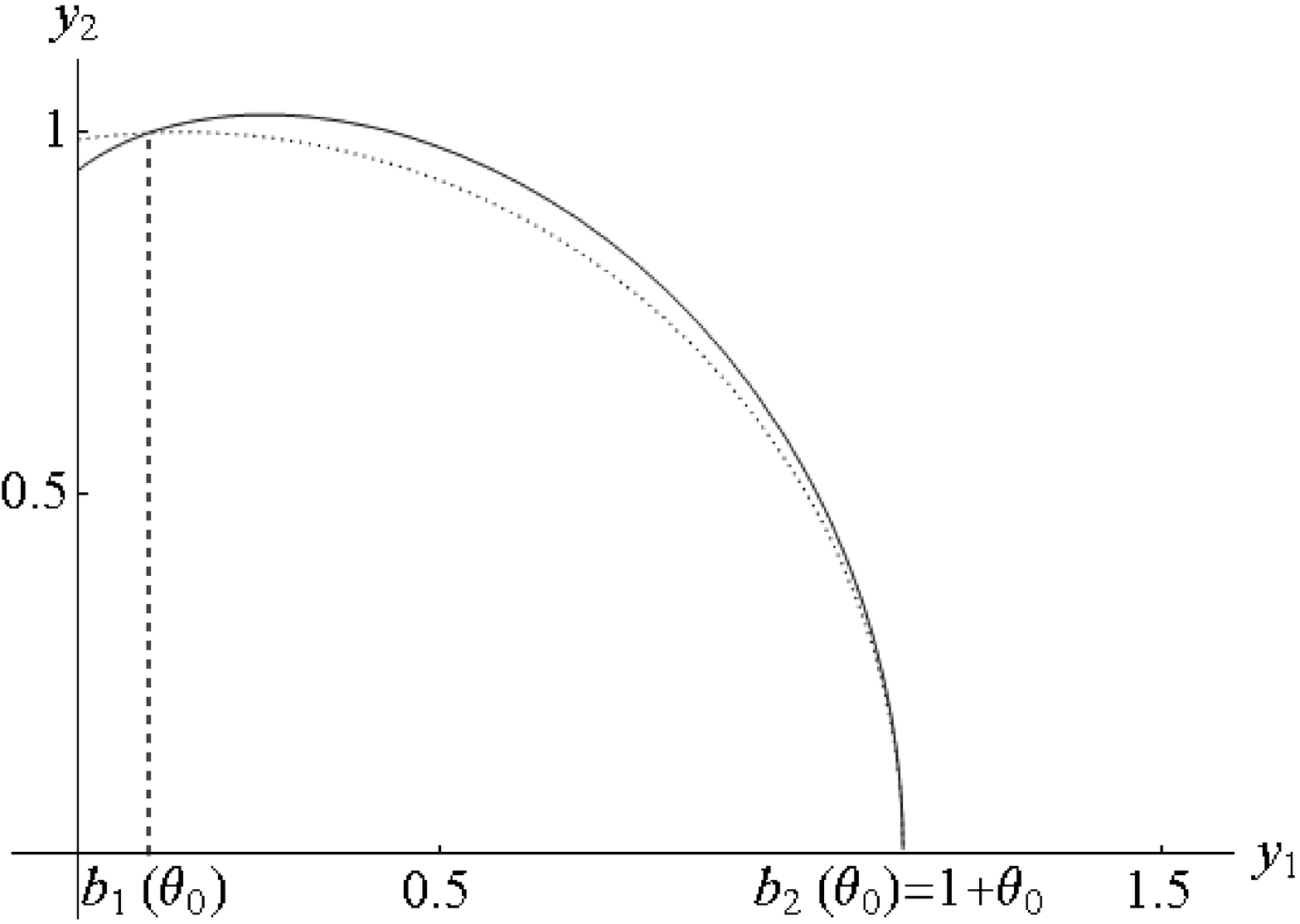}}
  \hbox to 9cm{\hss Figure 3:  $\theta=\theta_0$.\hss}
}}

\InsertBoxC{\box3}

\def\box4{\vbox{
  \epsfxsize=9cm
  \hbox{\epsfbox{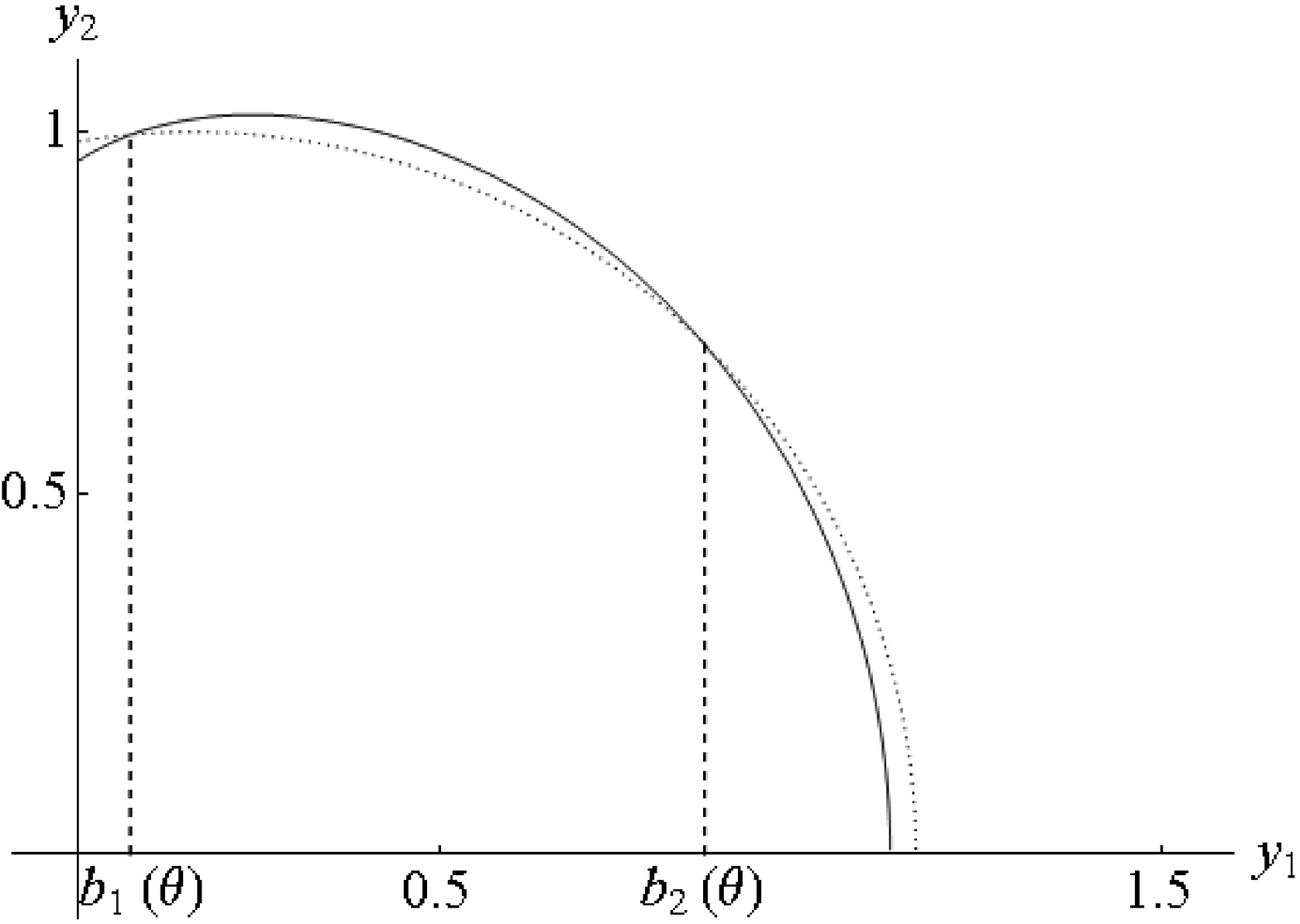}}
  \hbox to 9cm{\hss Figure 4:  $\theta_0<\theta<\theta_1$.\hss}
}}

\InsertBoxC{\box4}

\def\box5{\vbox{
  \epsfxsize=9cm
  \hbox{\epsfbox{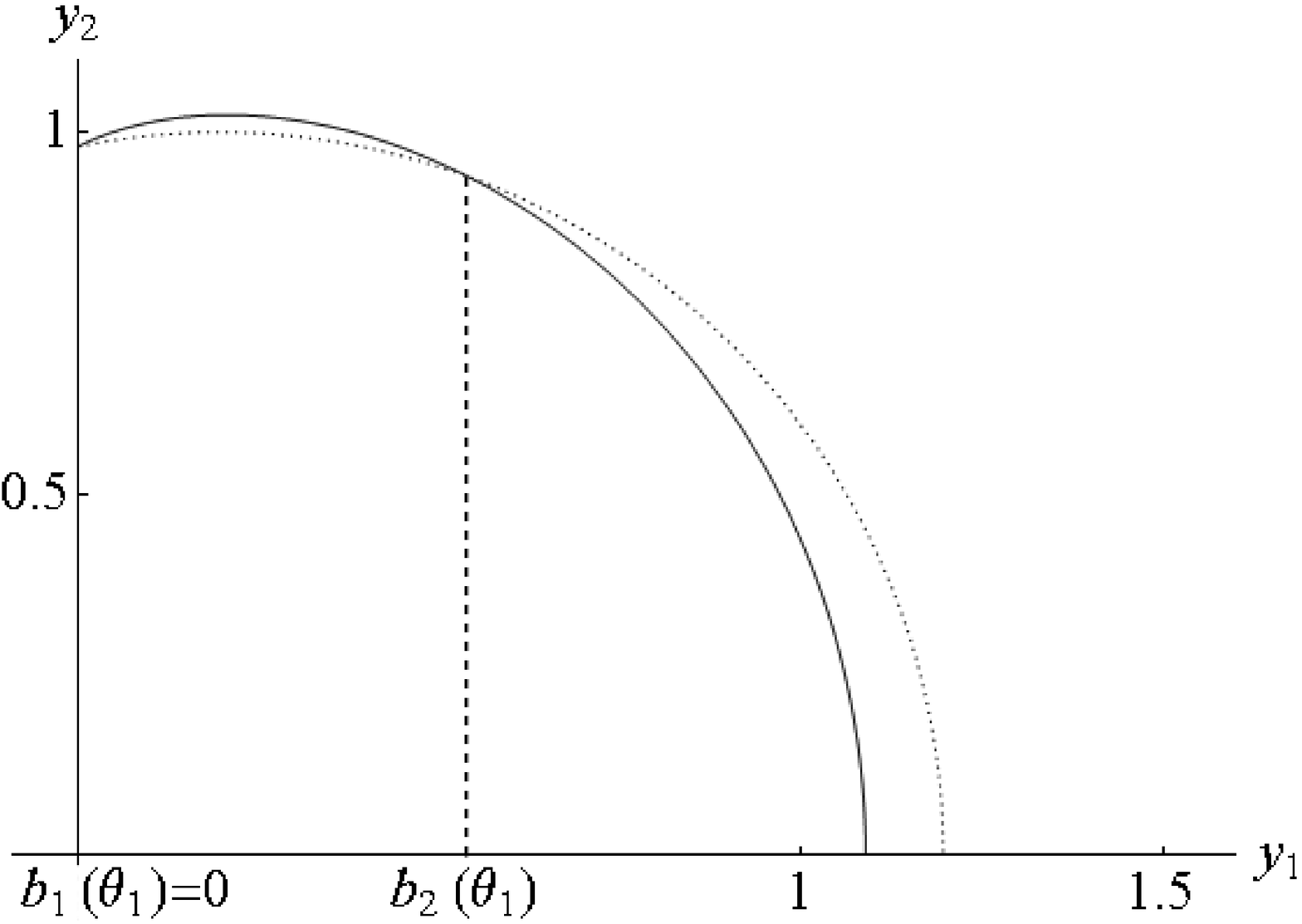}}
  \hbox to 9cm{\hss Figure 5:  $\theta=\theta_1$.\hss}
}}

\InsertBoxC{\box5}
\vskip3em
\noindent{\bf Acknowledgment.} The authors would like to thank Rick Laugesen for helpful discussions related to the computation  of the kernel $K(x,z)$, in the proof of  Theorem 5.\bigskip\eject
\centerline{\bf References}\bigskip
\item{[Ad]} Adams D.R. {\sl
A sharp inequality of J. Moser for higher order derivatives},
Ann. of Math. {\bf128} (1988), no. 2, 385--398. 
\smallskip
\item{[ABR]} Axler S., Bourdon P.,  Ramey W., {\sl Harmonic function theory, 2nd Ed.},  Springer-Verlag, New York,
2001.\smallskip
\item{[BFM]} Branson T.P., Fontana L., Morpurgo C., {\sl Moser-Trudinger and Beckner-Onofri's inequalities on the CR sphere}, 
(2007) submitted, arXiv:0712.3905.\smallskip
\item{[BM]} Brezis H., Merle F., {\sl Uniform estimates and blow-up behavior for solutions of in two dimensions}, Comm. Part. Diff. Eq. {\bf16} (1991), 1223-1253.\smallskip
\item{[CRT]}  Cassani D., Ruf B.,  Tarsi C., {\sl Best constants in a borderline case of second-order Moser type inequalities}, Ann. Inst. H. Poincaré Anal. Non Linéaire {\bf 27} (2010),  73-93.\smallskip 
\item{[CY]} Chang S.-Y.A., Yang  P.C., {\sl  Conformal deformation of metrics on $S^2$}, J. Differential Geom. {\bf 27} (1988), 259-296.
\smallskip\item{[Ci1]} Cianchi A., {\sl Moser-Trudinger trace inequalities}, Adv. Math. {\bf 217} (2008), 2005-2044.\smallskip
\item{[Ci2]} Cianchi A., {\sl Moser-Trudinger inequalities without boundary conditions and isoperimetric problems}, Indiana Univ. Math. J. {\bf 54} (2005), 669-705. 
\smallskip\item{[Fo]}  Fontana L.,  {\sl Sharp borderline Sobolev inequalities on compact Riemannian manifolds}, Comment. Math. Helv. 
{\bf68} (1993), 415-454.
\smallskip\item{[FM]} Fontana L., Morpurgo C., {\sl Adams inequalities on measure spaces}, (2009),  arXiv:0906.03, to appear in Adv. Math.
\smallskip\item{[GT]} Gilbarg D., Trudinger N.S., {\sl Elliptic Partial Differential Equations of Second Order}, 2nd ed., Springer-Verlag, New York, 1983.\smallskip
\item{[HP]} Harvey F.R., Polking J.C. {\sl The $\bar\partial$-Neumann solution to the inhomogeneous Cauchy-Riemann equation in the ball in $\C^n$}, Trans. Amer. Math. Soc. {\bf 281} (1984), 587-613.
\smallskip
\item{[KK]} Kang H., Koo H., {\sl Estimates of the Harmonic Bergman Kernel on Smooth Domains}, J. Funct. Anal. {\bf 185} (2001),  220-239. \smallskip
\item{[Ke]} Kershaw D., {\sl Some extensions of W. Gautschi's inequalities for the gamma function}, Math. Comp. {\bf 41} (1983),
607-611.\smallskip
\item{[Le]}  Leckband M., {\sl Moser's inequality on the ball $B^n$ for functions with mean value zero}, Comm. Pure Appl. Math. {\bf58} (2005), 789-798.\smallskip 
\item{[Li1]} Ligocka E., {\sl
On the reproducing kernel for harmonic functions and the space of Bloch harmonic functions on the unit ball in $\R^n$},  
Studia Math. {\bf 87} (1987), 23-32. \smallskip
 \item{[Li2]} Ligocka E., {\sl
Corrigendum to the paper: ``On the reproducing kernel for harmonic functions and the space of Bloch harmonic functions on the unit ball in $\R^n$''},
Studia Math. {\bf 101} (1992),  319. \smallskip
 \item{[Mo]} Moser J., {\sl A sharp form of an inequality by N. Trudinger}, Indiana Univ. Math. J. {\bf20} (1970/71), 1077-1092.\smallskip
\item{[PPCR]} Pankka P., Poggi-Corradini P., Rajala K. {\sl Sharp exponential integrability for traces of monotone Sobolev functions}, Nagoya Math. J. {\bf 192} (2008), 137-149. 
\bigskip

\noin Luigi Fontana \hskip19em Carlo Morpurgo

\noin Dipartimento di Matematica ed Applicazioni \hskip5.5em Department of Mathematics 
 
\noin Universit\'a di Milano-Bicocca\hskip 13em University of Missouri, Columbia

\noin Via Cozzi, 53 \hskip 19.3em Columbia, Missouri 65211

\noin 20125 Milano - Italy\hskip 16.6em USA 
\smallskip\noin luigi.fontana@unimib.it\hskip 15.3em morpurgoc@missouri.edu

\end